\documentclass[11pt,a4wide,reqno]{article}
\usepackage{verbatim}
\usepackage{graphics}
\usepackage[pdftex]{graphicx}
\usepackage{amsmath}
\usepackage{amssymb,latexsym}
\usepackage[only,fivrm]{rawfonts}
\usepackage{bbm}
\usepackage{enumerate}
\usepackage{amsthm}
\usepackage{multirow}
\usepackage{nccmath}
\usepackage{indentfirst}
\newcommand{\BR}{\Lambda_n}
\newcommand{\hb}{\hat{\beta}_n}
\newcommand{\GCM}{\widetilde{\Lambda}_n}
\newcommand{\ho}{\tilde{\lambda}_n}
\newcommand{\hoMLE}{\hat{\lambda}_n}

\newcommand{\fo}{\tilde{f}_{n}}

\def\inprob{\xrightarrow{p}}

\def\eqd{\,{\buildrel d \over =}\,}
\def\law{\xrightarrow{d}}
\renewcommand{\O}{\ensuremath{{\cal O}}}

\newtheorem{theorem}{Theorem}
\newtheorem{lemma}{Lemma}
\newtheorem{remark}{\emph{Remark}}
\newtheorem{corollary}{Corollary}

\usepackage{color}

\def\argmin{\mathop{\text{\sl\em argmin}}}
\def\argmax{\mathop{\text{\sl\em argmax}}}
\begin{document}
\setlength{\baselineskip}{18pt}

%\begin{comment}
\title{Shape constrained nonparametric estimators of the baseline distribution in Cox proportional hazards model}
\author {HENDRIK P.~LOPUHA\"{A} and GABRIELA F.~NANE\\[3mm]
\emph{Delft Institute of Applied Mathematics, Delft University of Technology}}
\date{\today}
\maketitle

\begin{abstract}
We investigate nonparametric estimation of a monotone baseline hazard and
a decreasing baseline density within the Cox model.
Two estimators of a nondecreasing baseline hazard function are proposed.
We derive the nonparametric maximum likelihood estimator and consider
a Grenander type estimator, defined as the left-hand slope of the greatest convex minorant of the Breslow estimator.
We demonstrate that the two estimators are strongly consistent and asymptotically equivalent and derive their common limit distribution at a fixed point.
Both estimators of a nonincreasing baseline hazard and their asymptotic properties are acquired in a similar manner.
Furthermore, we introduce a Grenander type estimator for a nonincreasing baseline density, defined
as the left-hand slope of the least concave majorant of an estimator of
the baseline cumulative distribution function, derived from the Breslow estimator.
We show that this estimator is strong consistent and derive its asymptotic distribution at a fixed point.

\bigskip


\bigskip

\noindent
\emph{Running headline:} Shape constrained estimation in the Cox model
\end{abstract}

\section{Introduction}
\label{sec:introduction}
\noindent
Shape constrained nonparametric estimation dates back to the $1950$s.
The milestone paper of Grenander~\cite{grenander:1956} introduced the maximum likelihood estimator of a nonincreasing density,
while Prakasa Rao~\cite{prakasarao:1969} derived its asymptotic distribution at a fixed point. Similarly, the maximum likelihood estimator of a monotone hazard function has been proposed by Marshall and Proschan~\cite{marshallproschan:1965}
and its asymptotic distribution was determined in~\cite{prakasarao:1970}.
Other estimators have been proposed and despite the high interest and applicability, the difficulty in the derivation of the asymptotics was a major drawback. Shape constrained estimation was revived by Groeneboom~\cite{groeneboom:1985},
who proposed an alternative for Prakasa Rao's bothersome type of proof.
Groeneboom's approach employs a so-called inverse process and makes use of the Hungarian embedding~\cite{komlosmajortusnady:1957}.
Once such an embedding is available, it enables the derivation of the asymptotic distribution of the considered estimator. This is the case, for example, when estimating a monotone density or hazard function from right-censored observations, as proposed by
Huang and Zhang~\cite{huangzhang:1994} and Huang and Wellner~\cite{huangwellner:1995}.
Another setting for deriving the asymptotic distribution, that does not require a Hungarian embedding, was later provided by the limit theorems in~\cite{kimpollard:1990}.
Their cube root asymptotics are based on a functional limit theorem for empirical processes.

The present paper treats
the estimation of a
monotone baseline hazard and a decreasing baseline density in the Cox model.
Ever since the model was introduced (see~\cite{cox:1972}) and
in particular,
since the asymptotic properties of the proposed estimators were first derived by Tsiatis~\cite{tsiatis:1981},
the Cox model is the classical survival analysis framework for incorporating covariates in the study of a lifetime distribution.
The hazard function is of particular interest in survival analysis, as it represents an important feature of the time course of a process under study, e.g., death or a certain disease. The main reason lies in its ease of interpretation and in the fact that the hazard function takes into account ageing, while, for example, the density function does not. Times to death, infection or development of a disease of interest in most survival analysis studies are observed to have a nondecreasing baseline hazard. Nevertheless, the survival time after a successful medical treatment is usually modeled using a nonincreasing hazard. An example of nonincreasing hazard is presented in Cook et al.~\cite{cooketal:1998},
where the authors concluded that the daily risk of pneumonia decreases with increasing duration of stay in the intensive care unit.

Chung and Chang~\cite{chungchang:1994} consider a maximum likelihood estimator of a nondecreasing baseline hazard function in the Cox model,
adopting the convention that each censoring time is equal to its preceding observed event time.
They prove consistency, but no distributional theory is available.
We consider a maximum likelihood estimator $\hoMLE$ of a monotone baseline hazard function, which imposes no extra assumption on the censoring times.
This estimator differs from the one in~\cite{chungchang:1994} and has a higher likelihood.
Furthermore, we introduce a Grenander type estimator for a monotone baseline hazard function based on the well-known baseline cumulative hazard estimator, the Breslow estimator $\BR$~\cite{cox:1972}.
The nondecreasing baseline hazard estimator $\ho$ is defined as the left-hand slope of the greatest convex minorant (GCM) of~$\BR$.
Similarly, a nonincreasing baseline estimator is characterized as the left-hand slope of the least concave majorant (LCM) of~$\BR$.
It is noteworthy that, just as in the no covariates case (see~\cite{huangwellner:1995}),
the two monotone estimators are different, but are shown to be asymptotically equivalent.
Additionally, we introduce a nonparametric estimator for a nonincreasing baseline density.
An estimator $F_n$ for the baseline distribution function is based on the Breslow estimator
and next, the baseline density estimator $\fo$ is defined as the left-hand slope of the LCM of $F_n$.
The treatment of the maximum likelihood estimator for a nonincreasing baseline density
is much more complex and is deferred to another paper.
For the remaining three estimators, we show that they converge at rate $n^{1/3}$ and we establish their limit distribution.
Since, to the authors best knowledge, there does not exist a Hungarian embedding for the Breslow estimator,
our results are based on the theory in~\cite{kimpollard:1990} and an argmax continuous mapping theorem in~\cite{huangwellner:1995}.

The paper is organized as follows.
In Section~\ref{sec:notation} we introduce the model and state our assumptions.
The formal characterization of the maximum likelihood estimator $\hoMLE$ is given in Lemmas~\ref{lemma:char incr}
and~\ref{lemma:char ML decreasing}.
Our main results concerning the asymptotic properties of the proposed estimators are gathered in Section~\ref{sec:main_results}.
Section~\ref{sec:consistency} is devoted to proving the strong consistency results of the paper. The strong uniform consistency of the Breslow estimator in Theorem~\ref{theorem:conv_breslow} and of the baseline cumulative distribution estimator $F_n$ in Corollary~\ref{cor:consistency_F_n}, emerge as necessary results.
These results are preceded by three preparatory lemmas, that establish properties of functionals in terms of which derivations thereof can be expressed.
In order to prepare the application of results from~\cite{kimpollard:1990},
in Section~\ref{sec:limit_distr_preliminary}
we introduce the inverses of the estimators in terms of minima and maxima
of random processes and obtain the limiting distribution of these processes.
Finally, in Section~\ref{sec:limit distribution}, we derive the asymptotic distribution of the estimators, at a fixed point.
The proofs of some preparatory lemmas are deferred to an appendix, which is available in the online Supporting Information.

\section{Definitions and assumptions}
\label{sec:notation}
\noindent
Let the observed data consist of independent identically distributed triplets~$\left(T_i,\Delta_i,Z_i\right)$,
with $i=1,2,\ldots,n$, where~$T_i$ denotes the follow-up time, with a corresponding censoring indicator~$\Delta_i$ and covariate vector $Z_i\in \mathbbm{R}^p$.
A generic follow-up time is defined by $T=\text{min}\left(X,C\right)$,
where $X$ represents the event time and~$C$ is the censoring time.
Accordingly, $\Delta=\left\{X\leq C\right\}$, where $\{\cdot\}$ denotes the indicator function.
The event time $X$ and censoring time $C$ are assumed to be
conditionally independent given $Z$, and the censoring mechanism is assumed to be non-informative.
The covariate vector $Z \in \mathbb{R}^p$ is assumed to be time invariant.

Within the Cox model, the distribution of the event time is related to the corresponding covariate by
\begin{equation}
\label{def:cox model}
\lambda\left(x|z\right)=\lambda_0(x)\,\text{e}^{\beta_0' z},
\end{equation}
where $\lambda\left(x|z\right)$ is the hazard function for an individual with covariate vector $z\in \mathbb{R}^p$,
$\lambda_0$ represents the baseline hazard function and $\beta_0\in \mathbbm{R}^p$
is the vector of the underlying regression coefficients.
Conditionally on $Z=z$,
the event time $X$ is assumed to be a nonnegative random variable with an absolutely continuous distribution function $F(x|z)$
with density $f(x|z)$.
The same assumptions hold for the censoring variable $C$ and its distribution function~$G$.
The distribution function of the follow-up time $T$ is denoted by $H$.
We will assume the following conditions, which are
commonly employed in deriving large sample properties of Cox proportional hazards estimators
(e.g., see~\cite{tsiatis:1981}).
\begin{description}
\item (A1)
Let $\tau_F,\tau_G$ and $\tau_H$ be the end points of the support of $F,G$ and $H$ respectively.
Then
\[
\tau_{H}=\tau_G<\tau_F\leq \infty.
\]
\item (A2)
There exists $\varepsilon>0$ such that
\[
\sup_{|\beta-\beta_0 |\leq \varepsilon} \mathbbm{E}\left[ |Z|^2\, \text{e}^{2\beta' Z}\right]<\infty,
\]
where $|\cdot|$ denotes the Euclidean norm.
\end{description}

\subsection{Increasing baseline hazard}
\label{subsec:increasing haz}
\noindent
Let $\Lambda(x|z)=-\log (1-F(x|z))$ be the cumulative hazard function.
Then, from~\eqref{def:cox model} it follows that $\Lambda(x|z)=\Lambda_0(x)\exp(\beta_0'z)$, where
$\Lambda_0(x)=\int_0^x \lambda_0(u)\,\mathrm{d}u$ denotes the
baseline cumulative hazard function.
When $G$ has a density $g$, then together with the relation $\lambda=f/(1-F)$,
the likelihood becomes
\[
\begin{split}
&
\prod_{i=1}^n
\Big[
f(T_i\mid Z_i)(1-G(T_i\mid Z_i))
\Big]^{\Delta_i}
\Big[
g(T_i\mid Z_i)(1-F(T_i\mid Z_i))
\Big]^{1-\Delta_i}\\
&=
\prod_{i=1}^n
\lambda(T_i\mid Z_i)^{\Delta_i}
\exp
\Big(-\Lambda(T_i\mid Z_i)\Big)
\times
\prod_{i=1}^n
\Big[1-G(T_i\mid Z_i)\Big]^{\Delta_i}
g(T_i\mid Z_i)^{1-\Delta_i}.
\end{split}
\]
The term with $g$ does not involve the baseline distribution and can be treated as a constant term.
Therefore, one essentially needs to maximize
\[
\prod_{i=1}^n
\lambda(T_i\mid Z_i)^{\Delta_i}
\exp\Big(-\Lambda(T_i\mid Z_i)\Big)
=
\prod_{i=1}^n
\Big[\lambda_0(T_i)\mathrm{e}^{\beta_0'Z_i}\Big]^{\Delta_i}
\exp\Big(-\mathrm{e}^{\beta_0'Z_i}\Lambda_0(T_i)\Big).
\]
This leads to the following (pseudo) loglikelihood,
written as a function of $\beta\in \mathbb{R}^p$ and $\lambda_0$,
\begin{equation}
\label{eq:loglikelihood}
\sum_{i=1}^n
\left[
\Delta_i\log\lambda_0(T_i) +\Delta_i\beta' Z_i-\text{e}^{\beta' Z_i}\Lambda_0(T_{i})
\right].
\end{equation}
\begin{remark}
\label{rem:discrete censoring}
It may be worthwhile to note that if the censoring distribution is discrete,
the likelihood of $(T,\Delta,Z)$ can still be written as
\[
[f(T\mid Z)(1-G(T\mid Z))]^{\Delta}
[g(T\mid Z)(1-F(T\mid Z))]^{1-\Delta},
\]
where $g(y|z)=P(C=y|Z=z)$, which will lead to the same expression as in~\eqref{eq:loglikelihood}.
However, as we will make use of other results in the literature that are established under the assumption
of an absolutely continuous censoring distribution (e.g., from~\cite{tsiatis:1981}), we do not further investigate the behavior of our estimators
in the case of a discrete censoring distribution.
\end{remark}
For $\beta\in \mathbb{R}^p$ fixed,
we first consider maximum likelihood estimation for a nondecreasing~$\lambda_0$.
This requires the maximization of~\eqref{eq:loglikelihood} over all nondecreasing $\lambda_0$.
Let $T_{(1)}<T_{(2)}<\cdots<T_{(n)}$ be the ordered follow-up times and,
for~$i=1,2,\ldots,n$, let $\Delta_{(i)}$ and~$Z_{(i)}$ be the censoring indicator and covariate vector corresponding to $T_{(i)}$.
The characterization of the maximizer~$\hoMLE(x;\beta)$ can be described by means of the processes
\begin{equation}
\label{def:W_n}
W_n(\beta,x)
=
\int \left( \text{e}^{\beta' z} \int_0^x \{ u\geq s \}\,\mathrm{d}s \right)\, \mathrm{d}\mathbbm{P}_n(u,\delta,z),
\end{equation}
and
\begin{equation}
\label{def:V_n}
V_n(x)
=
\int
\delta\{u< x\}\,\mathrm{d}\mathbbm{P}_n(u,\delta,z),
\end{equation}
with $\beta\in \mathbb{R}^p$ and $x\geq 0$,
where $\mathbbm{P}_n$ is the empirical measure of
the $(T_i,\Delta_i,Z_i)$ and is given by the following lemma.
\begin{lemma}
\label{lemma:char incr}
For a fixed $\beta\in \mathbbm{R}^p$, let $W_n$ and $V_n$ be defined in~\eqref{def:W_n} and~\eqref{def:V_n}.
Then, the NPMLE $\hoMLE(x;\beta)$ of a nondecreasing baseline hazard function $\lambda_0$ is of the form
\[
\hoMLE(x;\beta)
=
\begin{cases}
0 & x < T_{(1)},\\
\hat{\lambda}_i & T_{(i)}\leq x < T_{(i+1)}, \text{ for }i=1,2,\ldots,n-1,\\
\infty & x\geq T_{(n)},
\end{cases}
\]
where $\hat{\lambda}_i$ is the left derivative of the greatest convex minorant at the point $P_i$
of the cumulative sum diagram consisting of the points
\[
P_j=\Big( W_n(\beta,T_{(j+1)})-W_n(\beta,T_{(1)}), V_n(T_{(j+1)}) \Big),
\]
for $j=1,2,\ldots,n-1$ and $P_0=(0,0)$.
Furthermore,
\begin{equation}
\label{eq:maxmin increasing}
\hat{\lambda}_i
=
\max_{1\leq s\leq i}\,\, \min_{i\leq t\leq n-1}
\frac{\sum_{j=s}^t \Delta_{(j)}}{\sum_{j=s}^t \left(  T_{(j+1)}-T_{(j)} \right) \sum_{l=j+1}^n \text{e}^{\beta' Z_{(l)}}},
\end{equation}
for $i=1,2,\ldots,n-1$.
\end{lemma}
\begin{proof}
Similar to~\cite{marshallproschan:1965} and Section~7.4 in~\cite{robertsonetal},
since $\lambda_0(T_{(n)})$ can be chosen arbitrarily large, we first consider the maximization over nondecreasing~$\lambda_0$
bounded by some $M>0$.
When we increase the value of $\lambda_0$ on an interval $(T_{(i-1)},T_{(i)})$, the terms $\lambda_0(T_{(i)})$
in~\eqref{eq:loglikelihood} are not changed,
whereas terms with $\Lambda_0(T_{(i)})$ will decrease the loglikelihood.
Since $\lambda_0$ must be nondecreasing, we conclude that the solution is a
nondecreasing step function, that is zero for $x<T_{(1)}$,
constant on $[T_{(i)},T_{(i+1)})$, for $i=1,2,\ldots,n-1$, and equal to $M$, for $x\geq T_{(n)}$.
Consequently, for $\beta\in \mathbb{R}^p$ fixed, the (pseudo) loglikelihood reduces to
\begin{equation}
\label{eq:objective function}
\begin{split}
L_\beta(\lambda_0)
&=
\sum_{i=1}^{n-1}
\Delta_{(i)}\log\lambda_0(T_{(i)})
-
\sum_{i=2}^n
\text{e}^{\beta' Z_{(i)}}
\sum_{j=1}^{i-1}
\left( T_{(j+1)}-T_{(j)}\right)
\lambda_0(T_{(j)})\\
&=
\sum_{i=1}^{n-1}
\left[
\Delta_{(i)}\log\lambda_0(T_{(i)})
-
\lambda_0(T_{(i)})
\left( T_{(i+1)}-T_{(i)}\right)
\sum_{l=i+1}^n
\text{e}^{\beta' Z_{(l)}}
\right].
\end{split}
\end{equation}
Maximization over $0\leq \lambda_0(T_{(1)})\leq\cdots\leq \lambda_0(T_{(n-1)})\leq M$
will then have a solution $\hoMLE^M(x;\beta)$
and by letting $M\to\infty$, we obtain the NPMLE $\hoMLE(x;\beta)$ for $\lambda_0$.

First, notice that the loglikelihood function in \eqref{eq:objective function} can also be written as
\begin{equation}
\label{eq:isotonic loglik}
\begin{split}
\sum_{i=1}^{n-1}
\Big[
s_i\log\lambda_0(T_{(i)})-\lambda_0(T_{(i)})
\Big] w_i,
\end{split}
\end{equation}
where, for $i=1,2,\ldots,n-1$,
\[
w_i=\left(  T_{(i+1)}-T_{(i)} \right) \sum_{l=i+1}^n \text{e}^{\beta' Z_{(l)}},
\]
and
\[
s_i=\frac{\Delta_{(i)}}{\left(  T_{(i+1)}-T_{(i)} \right) \sum_{l=i+1}^n \text{e}^{\beta' Z_{(l)}}}.
\]
As mentioned above, we first maximize over nondecreasing $\lambda_0$ bounded by some~$M$.
Since~$M$ can be chosen arbitrarily large, the problem of maximizing~\eqref{eq:isotonic loglik}
over $0\leq \lambda_0(T_{(1)})\leq\cdots\leq \lambda_0(T_{(n-1)})\leq M$ can be identified
with the problem solved in Example 1.5.7 in~\cite{robertsonetal}.
The existence of $\hoMLE^M$ is therefore immediate and is given by
\[
\hoMLE^M(x;\beta)
=
\begin{cases}
0 & x < T_{(1)},\\
\hat{\lambda}_i & T_{(i)}\leq x < T_{(i+1)}, \text{ for }i=1,2,\ldots,n-1,\\
M & x\geq T_{(n)},
\end{cases}
\]
where, as a result of Theorems~1.5.1 and 1.2.1 in~\cite{robertsonetal},
the value $\hat{\lambda}_i$ is the left derivative at $P_i$ of the GCM of the cumulative sum diagram (CSD) consisting of the points
\[
P_i=\left(\frac{1}{n}\sum_{j=1}^iw_j,\frac{1}{n}\sum_{j=1}^iw_js_j\right), \qquad i=1,2,\ldots,n-1,
\]
and $P_0=(0,0)$.
It follows that
\[
\begin{split}
\frac{1}{n}\sum_{j=1}^i w_j
&=
\sum_{j=1}^i
\left(T_{(j+1)}-T_{(j)}\right)
\frac1n
\sum_{l=1}^n
\{T_l\geq T_{(j+1)}\}\text{e}^{\beta'Z_l}\\
&=
\int_{T_{(1)}}^{T_{(i+1)}}
\int
\{u\geq s\}\text{e}^{\beta'z}\,\mathrm{d}\mathbbm{P}_n(u,\delta,z)\,\mathrm{d}s
=
W_n(\beta,T_{(i+1)})-W_n(\beta,T_{(1)}).
\end{split}
\]
For the $y$-coordinate of the CSD, notice that
\[
\frac{1}{n}\sum_{j=1}^i w_js_j
=
\frac{1}{n} \sum_{j=1}^i \Delta_{(j)}
=
\frac{1}{n}\sum_{j=1}^n \{T_j\leq T_{(i)},\Delta_j=1\}
=
V_n(T_{(i+1)}).
\]
By letting $M\to\infty$, we obtain the NPMLE $\hoMLE(\beta,x)$ for $\lambda_0$.
The max-min formula in \eqref{eq:maxmin increasing} follows from Theorem~1.4.4 in~\cite{robertsonetal}.
\end{proof}

\noindent
\emph{Remark:}
From the characterization given in Lemma~\ref{lemma:char incr}, it can be seen that the GCM of the CSD  only changes slope
at points corresponding to uncensored observations, which means that~$\hoMLE(x;\beta)$ is constant between successive uncensored follow-up times.
Moreover, similar to the reasoning in the proof of Lemma~\ref{lemma:char incr}, it follows that $\hoMLE(x;\beta)$ maximizes~\eqref{eq:loglikelihood}.
The reason to provide the characterization in Lemma~\ref{lemma:char incr} in terms of all follow-up times
is that this facilitates the treatment of the asymptotics for this estimator.
Finally, for the solution~$\hoMLE^M(x;\beta)$, on the interval $[T_{(n)},\tau_H)$, in principle one could take any value
between~$\hat{\lambda}_{n-1}$ and~$M$.
This means that for $\hoMLE(x;\beta)$, on the interval $[T_{(n)},\tau_H)$, one could take any value larger than
$\hat{\lambda}_{n-1}$.

\bigskip

In practice, one also has to estimate $\beta_0$.
The standard choice is $\hb$, the maximizer of the partial likelihood function
\[
\prod_{l=1}^m\frac{\text{e}^{\beta'Z_j}}{\sum_{j=1}^n\{T_j\geq X_{(i)}\}\text{e}^{\beta' Z_j}},
\]
as proposed by Cox~\cite{cox:1972,cox:1975}, where $X_{(1)}<X_{(2)}<\cdots<X_{(m)}$ denote the ordered, observed event times.
Since the maximum partial likelihood estimator~$\hb$ for $\beta_0$ is asymptotically efficient under mild conditions
and because the amount of information on $\beta_0$ lost through lack
of knowledge of $\lambda_0$ is usually small
(see e.g.,\cite{efron:1977,oakes:1977,slud:1982}), we do not pursue joint maximization of~\eqref{eq:loglikelihood} over
nondecreasing $\lambda_0$ and $\beta_0$.
We simply replace $\beta$ in $\hoMLE(x;\beta)$ by $\hb$, and we propose $\hoMLE(x)=\hoMLE(x;\hb)$ as
our estimator for $\lambda_0$.

Note that $\hoMLE$ is different from the estimator derived in~\cite{chungchang:1994},
where each censoring time is taken equal to the preceding observed event time.
This leads to a CSD that is slightly different from the one in Lemma~\ref{lemma:char incr}.
However, it can be shown that both estimators have the same asymptotic behavior.
Furthermore, if we take all covariates equal to zero,
the model coincides with the ordinary random censorship model with a nondecreasing hazard function as considered in~\cite{huangwellner:1995}.
The characterization in Lemma~\ref{lemma:char incr}, with all $Z_l\equiv0$, differs slightly from the one in Theorem~3.2 in~\cite{huangwellner:1995}.
Their estimator seems to be the result of maximization of~\eqref{eq:loglikelihood} over left-continuous $\lambda_0$
that are constant between follow-up times.
Although this estimator does not maximize~\eqref{eq:loglikelihood} over all nondecreasing $\lambda_0$,
the asymptotic distribution will turn out to be the same as that of $\hoMLE$, for the special case of no covariates.
The computation of joint maximum likelihood estimates for $\beta$ and $\lambda_0$ is considered in~\cite{huijankowski:2012},
who also developed an R package to compute the estimates.

\begin{figure}[htp]
 \includegraphics[scale=0.4]{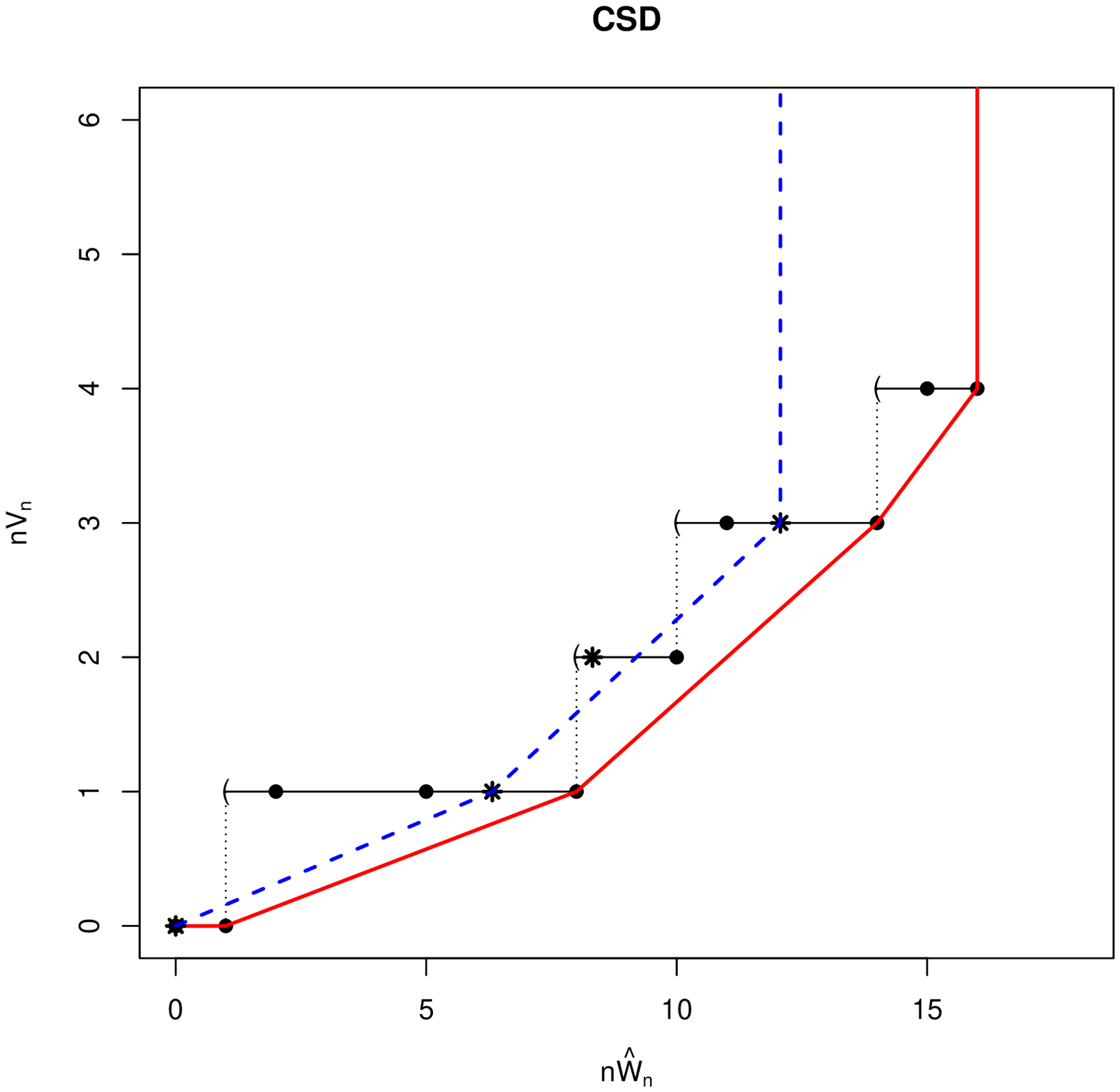}
 \includegraphics[scale=0.4]{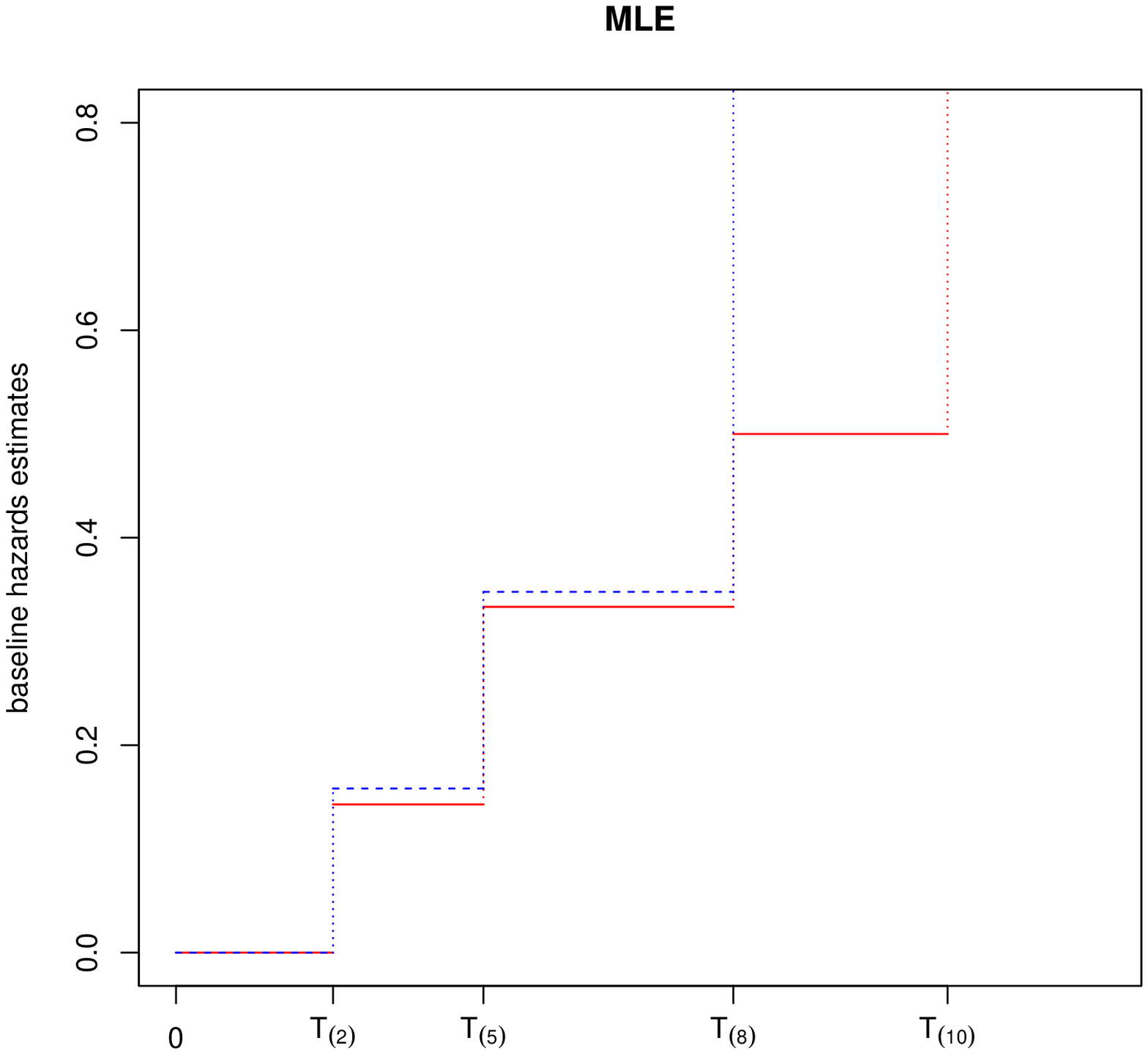}
\caption{The cumulative sum diagrams along with their GCM (left panel)
and the corresponding estimates of a nondecreasing baseline hazard (right panel).
Black points and solid curves correspond to the estimator in Lemma~\ref{lemma:char incr};
star points and dashed curves correspond to the estimator in~\cite{chungchang:1994}.}
 \label{fig:NPMLE}
\end{figure}

To illustrate the computation of the estimator described in Lemma~\ref{lemma:char incr},
consider an artificial survival dataset consisting of $10$ follow-up times,
with only $T_{(2)},T_{(5)},T_{(6)}$, and $T_{(8)}$ being observed event times.
In Figure~\ref{fig:NPMLE} we illustrate the construction of the proposed estimator and compare the resulting estimate
with the one suggested in~\cite{chungchang:1994}.
In order to compare the CSD of both estimates, the coordinates of the CSD described in
Lemma~\ref{lemma:char incr} have been multiplied with a factor $n$, which obviously leads to the same slopes.
Figure~\ref{fig:NPMLE} displays the points of the CSD (black points) and the GCM (solid curve) in the left panel.
The horizontal segments are generated by $(nW_n(\hb,x)-nW_n(\hb,T_{(1)}),nV_n(x))$ for $x\geq T_{(1)}$.
Note that the process~$nV_n$ has a jump of size 1 right after a point $P_j$ that corresponds to an observed event time.
Taking left derivatives then yield jumps of $\hoMLE$ only at observed event times.
The right panel of Figure~\ref{fig:NPMLE} displays the corresponding graph of $\hoMLE$ (solid curve).
The jumps of~$\hoMLE$ in the right panel correspond to the changes of slope of the GCM at the points
$P_1,P_4$ and $P_7$ in the left panel and occur at the event times $T_{(2)},T_{(5)}$, and $T_{(8)}$.
The height of the horizontal segments in the right panel corresponds to the slopes of the GCM in the left panel.
For comparison we have added the CSD (star points) and the corresponding GCM (dashed curve) of the estimator derived in~\cite{chungchang:1994}
in the left panel and the resulting estimator in the right panel (dashed curve).
Note that shifting the censoring times back to the nearest previous event time, as suggested in~\cite{chungchang:1994},
pushes points in the CSD, that correspond to event times, to the left.
As a consequence this yields steeper slopes in the left panel and hence a larger estimate of the hazard in the right panel.

\bigskip

\[
\text{Figure~\ref{fig:NPMLE} about here.}
\]

\bigskip

Another possibility to estimate a nondecreasing hazard is to construct a Grenander type estimator,
i.e., consider an unconstrained estimator $\Lambda_n$ for the cumulative hazard $\Lambda_0$
and take the left derivative of the GCM as an estimator of $\lambda_0$.
Several isotonic estimators are of this form (see e.g.,~\cite{grenander:1956, brunk:1958, huangwellner:1995, durot:2007}).
Breslow~\cite{cox:1972} proposed
\begin{equation}
\label{def:breslow}
\BR(x)=
\sum_{i|X_{(i)}\leq x}\dfrac{d_i}{\sum_{j=1}^n\{T_j\geq X_{(i)}\}\text{e}^{\hb' Z_j}},
\end{equation}
as an estimator for $\Lambda_0$,
where $d_i$ is the number of events at $X_{(i)}$
and $\hb$ is the maximum partial likelihood estimator of the regression coefficients.
The estimator $\BR$ is most commonly referred to as the Breslow estimator.
In the case of no covariates, i.e., $\beta=0$, the NPMLE estimate
of an increasing hazard rate has been illustrated in~\cite{huangwellner:1995}.

Following the derivations in~\cite{tsiatis:1981}, it can be inferred that
\begin{equation}
\label{eq:relation base haz}
\lambda_0(x)
=
\frac{ \mathrm{d}H^{uc}(x)/\mathrm{d}x }{\mathbb{E}\left[\{T\geq x\}\exp(\beta_0'Z)\right]},
\end{equation}
where $H^{uc}(x)=\mathbbm{P}(T\leq x, \Delta=1)$ is the sub-distribution function of the uncensored observations.
Consequently, it can be derived that
\begin{equation}
\label{eq:relation cum base haz}
\Lambda_0(x)=\int \frac{\delta\{u\leq x\}}{\mathbb{E}\left[\{T\geq x\}\exp(\beta_0'Z)\right]}\,\mathrm{d}P(u,\delta,z),
\end{equation}
where $P$ is the underlying probability measure corresponding to the distribution of $(T,\Delta,Z)$.
From (A1), it follows that $\Lambda_0(\tau_H)<\infty$.
In view of the above expression, an intuitive baseline cumulative hazard estimator is obtained by
replacing the expectations in~\eqref{eq:relation cum base haz} by averages and by plugging in~$\hb$,
which yields exactly the Breslow estimator in~\eqref{def:breslow}.
As a Grenander type estimator for a nondecreasing hazard,
we propose the left derivative $\ho$ of the greatest convex minorant $\GCM$ of $\BR$.
This estimator is different from $\hoMLE$ for finite samples, but we will show that
both estimators are asymptotically equivalent.
For the special case of no covariates, this coincides with the results in~\cite{huangwellner:1995}.
\subsection{Decreasing baseline hazard}
\label{subsec:decreasing haz}
\noindent
A completely similar characterization is provided for the NPMLE of a nonincreasing baseline hazard function.
As in the nondecreasing case, one can argue that the loglikelihood is maximized by
a decreasing step function that is constant on $(T_{(i-1)},T_{(i)}]$, for $i=1,2,\ldots,n$, where $T_{(0)}=0$.
In this case, the loglikelihood reduces to
\[
L_\beta(\lambda_0)
=
\sum_{i=1}^{n}
\left[
\Delta_{(i)}\log\lambda_0(T_{(i)})
-
\lambda_0(T_{(i)})
\left( T_{(i)}-T_{(i-1)}\right)
\sum_{l=i}^n
\text{e}^{\beta' Z_{(l)}}
\right],
\]
which is maximized over all $\lambda_0(T_{(1)})\geq \cdots \geq \lambda_0(T_{(n)})\geq 0$.
The solution is characterized by the following lemma.
The proof of this lemma is completely similar to that of Lemma~\ref{lemma:char incr}.
\begin{lemma}
\label{lemma:char ML decreasing}
For a fixed $\beta\in \mathbbm{R}^p$, let $W_n$ be defined in~\eqref{def:W_n} and
let
\begin{equation}
\label{def:Y_n}
Y_n(x)=\int\delta\{u\leq x\}\,\mathrm{d}\mathbb{P}_n(u,\delta,z).
\end{equation}
Then the NPMLE $\hoMLE(x;\beta)$ of a nonincreasing baseline hazard function $\lambda_0$ is given by
\[
\hoMLE(x;\beta)
=
\hat{\lambda}_i
\qquad
\text{for } x\in (T_{(i-1)},T_{(i)}],
\]
for $i=1,2,\ldots,n$, where $\hat{\lambda}_i$ is the left derivative of the least concave majorant (LCM) at the point $P_i$
of the cumulative sum diagram consisting of the points
\[
P_j=\Big( W_n(\beta,T_{(j)}), Y_n(T_{(j)}) \Big),
\]
for $j=1,2,\ldots,n$ and $P_0=(0,0)$.
Furthermore,
\[
\hat{\lambda}_i
=
\max_{1\leq s\leq i}\,\, \min_{i\leq t\leq n}
\frac{\sum_{j=s}^t \Delta_{(j)}}{\sum_{j=s}^t \left(  T_{(j)}-T_{(j-1)} \right) \sum_{l=j}^n \text{e}^{\beta' Z_{(l)}}},
\]
for $i=1,2,\ldots,n$.
\end{lemma}
Analogous to the nondecreasing case, for $x\geq T_{(n)}$, one can choose for $\hoMLE(x;\beta)$ any value smaller than $\hat{\lambda}_{n}$.
As before, we propose $\hoMLE(x)=\hoMLE(x;\hb)$ as an estimator for $\lambda_0$,
where $\hb$ denotes the maximum partial likelihood estimator for $\beta_0$.
Similar to the nondecreasing case, the Grenander type estimator $\ho$ for a nonincreasing $\lambda_0$ is defined as
the left-hand slope of the LCM of the Breslow estimator $\BR$, defined in~\eqref{def:breslow}.

An illustration of
the NPMLE of a decreasing baseline hazard function can be found in~\cite{nanetal:2012}, who investigated the hazard of patients with acute coronary syndrome.
Previous clinical trials indicated a decreasing risk pattern, which the authors confirmed by a test based on a bootstrap procedure.
The above estimate has been computed for 1200 patients undergoing early or selective invasive strategies, that were monitored for five years,
and their performance was evaluated by means of a simulation experiment.
The R code is available in the online version of their paper.

\subsection{Decreasing baseline density}
\label{subsec:mondens}
\noindent
Suppose one is interested in estimating a nonincreasing baseline density $f_0(\cdot)=f(\cdot|z=0)$.
One might argue that this problem is of less interest,
because the monotonicity assumption assumed for $z=0$ may no longer hold if one transforms the covariates by $a+bz$,
whereas the Cox model essentially remains unchanged.
Whereas the estimator for the baseline hazard remains monotone under such transformations,
this may no longer hold for the estimator of the baseline density.
Despite this drawback, we feel that the estimation of a nonincreasing baseline density may be of interest.

In this case, the corresponding baseline distribution function $F_0$ is concave
and it relates to the baseline cumulative hazard function $\Lambda_0$ as follows
\begin{equation}
\label{eq:relation F0}
F_0(x)=1-\text{e}^{-\Lambda_0(x)}.
\end{equation}
Hence, a natural estimator of the baseline distribution function is
\begin{equation}
\label{eq: relation Fn}
F_n(x)=1-\text{e}^{-\BR(x)},
\end{equation}
where $\BR$ is the Breslow estimator, defined in~\eqref{def:breslow}.
A Grenander type estimator~$\fo$ of a nonincreasing baseline density is defined as the left-hand slope of the LCM of $F_n$.
Recall that~$\BR$ depends on $\hb$ and $Z_1,Z_2,\ldots,Z_n$, and therefore the same holds for $F_n$ and $\fo$.

The derivation of the NPMLE for $f_0$ is much more complex than the previous estimators and its treatment is postponed to
a future manuscript.
In the special case of no covariates, the NPMLE $\hat f_n$ has first been derived in~\cite{huangzhang:1994}.
In~\cite{huangwellner:1995} a different characterization has been provided for $\hat f_n$  in terms of a self-induced cusum diagram
and it was shown that $\hat f_n$ and $\fo$ are asymptotically equivalent.

\section{Main results}
\label{sec:main_results}
\noindent
In this section, we state our main results.
The proofs are postponed to subsequent sections.
The next theorem provides pointwise consistency of the proposed estimators
at a fixed point~$x_0$ in the interior of the support.
Note that the results below imply that if~$x_0$ is a point of continuity of $\lambda_0$,
then $\hoMLE(x_0)\to \lambda_0(x_0)$ with probability one, and likewise for the other estimators.
\begin{theorem}
\label{theorem:consistency_estimators}
Assume that (A1) and (A2) hold.
\begin{enumerate}[(i)]
\item
Suppose that $\lambda_0$ is nondecreasing on $[0,\infty)$
and let $\hoMLE$ and $\ho$ be the estimators defined in Section~\ref{subsec:increasing haz}.
Then, for any~$x_0\in(0,\tau_H)$,
\[\begin{split}
\lambda_0(x_0-)\leq \liminf_{n \rightarrow \infty} \ \hat{\lambda}_n(x_0) \leq \limsup_{n \rightarrow \infty}
\hat{\lambda}_n(x_0) \leq \lambda_0(x_0+),\\
\lambda_0(x_0-)\leq \liminf_{n \rightarrow \infty} \ \tilde{\lambda}_n(x_0) \leq \limsup_{n \rightarrow \infty}
\tilde{\lambda}_n(x_0) \leq \lambda_0(x_0+),
\end{split}\]
with probability one,
where the values $\lambda_0(x_0-)$ and $\lambda_0(x_0+)$ denote the left and right limit at $x_0$.
\item
Suppose that $\lambda_0$ is nonincreasing on $[0,\infty)$
and let $\hoMLE$ and $\ho$ be the estimators defined in Section~\ref{subsec:decreasing haz}.
Then, for any~$x_0\in(0,\tau_H)$,
\[
\begin{split}
\lambda_0(x_0+)\leq \liminf_{n \rightarrow \infty} \hat{\lambda}_n(x_0)
&\leq
\limsup_{n \rightarrow \infty}\hat{\lambda}_n(x_0) \leq \lambda_0(x_0-),\\
\lambda_0(x_0+)\leq \liminf_{n \rightarrow \infty} \tilde{\lambda}_n(x_0)
&\leq
\limsup_{n \rightarrow \infty}\tilde{\lambda}_n(x_0) \leq \lambda_0(x_0-),
\end{split}
\]
with probability one.
\item
Suppose that $f_0$ is nonincreasing on $[0,\infty)$ and let $\fo$ be the estimator defined in Section~\ref{subsec:mondens}.
Then, for any~$x_0\in(0,\tau_H)$,
\[
f_0(x_0+)\leq \liminf_{n \rightarrow \infty} \tilde{f}_n(x_0) \leq \limsup_{n \rightarrow \infty}
\tilde{f}_n(x_0) \leq f_0(x_0-),
\]
with probability one,
where $f_0(x_0-)$ and $f_0(x_0+)$ denote the left and right limit at $x_0$.
\end{enumerate}
\end{theorem}
The following two theorems yield the asymptotic distribution of the monotone constrained baseline hazard estimators.
In order to keep notations compact, it becomes useful to introduce
\begin{equation}
\label{def:Phi}
\Phi(\beta,x)
=
\int \{u\geq x\}\, \text{e}^{\beta' z}\, \mathrm{d}P(u,\delta,z),
\end{equation}
for $\beta\in \mathbb{R}^p$ and $x\in \mathbb{R}$,
where $P$ is the underlying probability measure corresponding to the distribution of $(T,\Delta,Z)$.
Furthermore, by the argmin function we mean the supremum of times at which the minimum is attained.
Note that the limiting distribution and the rate of convergence coincide with the results commonly obtained
for isotonic estimators and differ from the corresponding quantities in the traditional central limit theorem.
The limiting distribution, usually referred to as the Chernoff distribution, has been tabulated
in~\cite{groeneboomwellner:2001}.
\begin{theorem}
\label{theorem:asymp_hazard_nondecr}
Assume (A1) and (A2) and let $x_0\in(0,\tau_H)$.
Suppose that~$\lambda_0$ is nondecreasing on $[0,\infty)$ and continuously differentiable in a neighborhood of $x_0$,
with $\lambda_0(x_0)\neq 0$ and $\lambda_0'(x_0)>0$.
Moreover, suppose that $H^{uc}(x)$ and $x\mapsto \Phi(\beta_0,x)$ are continuously differentiable in a neighborhood of $x_0$,
where $H^{uc}$ is defined below~\eqref{eq:relation base haz} and $\Phi$ is defined in~\eqref{def:Phi}.
Let~$\hoMLE$ and~$\ho$ be the estimators defined in Section~\ref{subsec:increasing haz}.
Then,
\begin{equation}
\label{eq:limit distribution increasing haz}
n^{1/3}
\left( \frac{\Phi(\beta_0,x_0)}{4\lambda_0(x_0)\lambda_0'(x_0)}\right)^{1/3}
\left(
\hoMLE(x_0)-\lambda_0(x_0)
\right)
\law
\underset{t\in \mathbbm{R}}{\operatorname{argmin}}\{\mathbbm{W}(t)+t^2\},
\end{equation}
where $\mathbbm{W}$ is standard two-sided Brownian motion originating from zero.
Furthermore,
\begin{equation}
\label{eq:asymp equiv MLE-LS}
n^{1/3}
\left(\ho(x_0)-\hoMLE(x_0) \right)
\inprob 0,
\end{equation}
so that the convergence in~\eqref{eq:limit distribution increasing haz} also holds with $\hoMLE$ replaced by $\ho$.
\end{theorem}
Let $\overline{\lambda}_n$ be the estimator considered in~\cite{chungchang:1994}, which has been proven to be consistent.
Completely similar to the proof of Theorem~\ref{theorem:asymp_hazard_nondecr} it can be shown that
\[
n^{1/3}
\left( \overline{\lambda}_n(x_0)-\hoMLE(x_0)
\right)
\inprob 0,
\]
so that the convergence in~\eqref{eq:limit distribution increasing haz} also holds with $\hoMLE$ replaced by $\overline{\lambda}_n$.
The next theorem establishes the same results as in Theorem~\ref{theorem:asymp_hazard_nondecr}, for the nonincreasing case.
\begin{theorem}
\label{theorem:asymp_hazard_nonincr}
Assume (A1) and (A2) and let $x_0\in(0,\tau_H)$.
Suppose that~$\lambda_0$ is nonincreasing on $[0,\infty)$ and continuously differentiable in a neighborhood of $x_0$,
with $\lambda_0(x_0)\neq 0$ and $\lambda_0'(x_0)<0$.
Moreover, suppose that $H^{uc}(x)$ and $x\mapsto \Phi(\beta_0,x)$ are continuously differentiable in a neighborhood of $x_0$,
where $H^{uc}$ is defined below~\eqref{eq:relation base haz} and $\Phi$ is defined in~\eqref{def:Phi}.
Let~$\hoMLE$ and $\ho$ be the estimators defined in Section~\ref{subsec:decreasing haz}.
Then,
\begin{equation}
\label{eq:limit distribution decreasing haz}
n^{1/3}
\left| \frac{\Phi(\beta_0,x_0)}{4\lambda_0(x_0)\lambda_0'(x_0)} \right|^{1/3}
\left(
\hoMLE(x_0)-\lambda_0(x_0)
\right)
\law
\underset{t\in \mathbbm{R}}{\operatorname{argmin}}\{\mathbbm{W}(t)+t^2\},
\end{equation}
where $\mathbbm{W}$ is standard two-sided Brownian motion originating from zero.
Furthermore,
\[
n^{1/3}\left( \ho(x_0)-\hoMLE(x_0) \right)\inprob 0,
\]
so that the convergence in~\eqref{eq:limit distribution decreasing haz} also holds with $\hoMLE$ replaced by $\ho$.
\end{theorem}
In the special case of no covariates, i.e., $Z\equiv0$, it follows that
$\Phi(\beta_0,x_0)=1-H(x_0)$, so that with the above results we recover Theorems~2.2 and~2.3 in~\cite{huangwellner:1995}.
If, in addition, one specializes to the case of no censoring, i.e., $\Phi(\beta_0,x_0)=1-H(x_0)=1-F(x_0)$, we recover Theorems~6.1 and~7.1 in~\cite{prakasarao:1970}.
The asymptotic distribution of the baseline density estimator is provided by the next theorem.
\begin{theorem}
\label{theorem:asymp_density}
Assume (A1) and (A2) and let $x_0\in(0,\tau_H)$.
Suppose that~$f_0$ is nonincreasing on $[0,\infty)$ and continuously differentiable in a neighborhood of $x_0$,
with $f_0(x_0)\neq 0$ and $f_0'(x_0)<0$.
Let $F_0$ be the baseline distribution function and
suppose that $H^{uc}(x)$ and $x\mapsto \Phi(\beta_0,x)$ are continuously differentiable in a neighborhood of $x_0$,
where $H^{uc}$ is defined below~\eqref{eq:relation base haz} and $\Phi$ is defined in~\eqref{def:Phi}.
Let $\fo$ be the estimator defined in Section~\ref{subsec:mondens}.
Then,
\[
n^{1/3}
\left| \frac{\Phi(\beta_0,x_0)}{4f_0(x_0)f_0'(x_0)[1-F_0(x_0)]} \right|^{1/3}
\left( \fo(x_0)-f_0(x_0)\right)
\law
\underset{t\in \mathbbm{R}}{\operatorname{argmin}}\{\mathbbm{W}(t)+t^2\},
 \]
where $\mathbbm{W}$ is standard two-sided Brownian motion originating from zero.
\end{theorem}
In the special case of no covariates, it follows that
\[
\frac{\Phi(\beta_0,x_0)}{1-F_0(x_0)}
=
\frac{1-H(x_0)}{1-F(x_0)}
=
1-G(x_0),
\]
so that the above result recovers Theorem~2.1 in~\cite{huangwellner:1995}.
If, in addition, one specializes to the case of no censoring, i.e., $G(x_0)=0$,
we recover Theorem~6.3 in~\cite{prakasarao:1969} and the corresponding result in~\cite{groeneboom:1985}.

\section{Consistency}
\label{sec:consistency}
\noindent
The strong pointwise consistency of the proposed estimators will be proven
using arguments similar to those in~\cite{robertsonetal} and~\cite{huangwellner:1995}.
First, define
\begin{equation}
\label{def:Phi_n}
\Phi_n(\beta,x)=\int \{u\geq x\}\, \text{e}^{\beta' z}\, \mathrm{d}\mathbbm{P}_n(u,\delta,z),
\end{equation}
for $\beta\in \mathbb{R}^p$ and $x\geq 0$
and note that the Breslow estimator in~\eqref{def:breslow} can also be represented as
\begin{equation}
\label{eq:relation Breslow Phi_n}
\Lambda_n(x)=\int \frac{\delta\{u\leq x\}}{\Phi_n(\hb,u)}\ \mathrm{d}\mathbbm{P}_n(u,\delta,z),
\qquad
x\geq 0.
\end{equation}
To establish consistency of the estimators, we first obtain some properties of
$\Phi_n$ and $\Phi$, as defined in~\eqref{def:Phi_n} and~\eqref{def:Phi}
and their first and second partial derivatives,
which by the dominated convergence theorem and conditions (A1) and~(A2) are given by
\[
\begin{split}
D^{(1)}(\beta,x)
&=
\frac{\partial\Phi(\beta,x)}{\partial\beta}
=
\int\{ u\geq x \}\, z \, \text{e}^{\beta' z}\, \mathrm{d} P(u,\delta,z)\,\in \mathbbm{R}^p,\\
D_n^{(1)}(\beta,x)
&=
\frac{\partial\Phi_n(\beta,x)}{\partial\beta}
=
\int\{ u\geq x \}\, z \, \text{e}^{\beta' z}\, \mathrm{d}\mathbbm{P}_n(u,\delta,z)\,\in \mathbbm{R}^p,\\
D^{(2)}(\beta,x)
&=
\frac{\partial^2\Phi(\beta,x)}{\partial\beta^2}
=
\int\{ u\geq x \}\, zz' \, \text{e}^{\beta' z}\, \mathrm{d}P(u,\delta,z)\,\in \mathbbm{R}^p\times\mathbbm{R}^p,\\
D_n^{(2)}(\beta,x)
&=
\frac{\partial^2\Phi_n(\beta,x)}{\partial\beta^2}
=
\int\{ u\geq x \}\, zz' \, \text{e}^{\beta' z}\, \mathrm{d}\mathbbm{P}_n(u,\delta,z)\,\in \mathbbm{R}^p\times\mathbbm{R}^p.
\end{split}
\]
In order to prove consistency, we need uniform bounds on $\Phi$ and its derivatives.
These are provided by the next lemma.
\begin{lemma}
\label{lemma:bounds Phin}
Suppose that (A2) holds for some $\varepsilon>0$. Then, for any $0<M<\tau_H$,
\begin{enumerate}[(i)]
\item
\[
0<
\inf_{x\leq M}
\inf_{|\beta-\beta_0|\leq\varepsilon} |\Phi(\beta,x)|
\leq
\sup_{x\in\mathbb{R}} \sup_{|\beta-\beta_0|\leq\varepsilon} |\Phi(\beta,x)|<\infty.
\]
\item
For any sequence $\beta_n^*$, such that $\beta^*_n \to \beta_0$ almost surely,
\[
0<\liminf_{n\to \infty}
\inf_{x\leq M}
|\Phi_n(\beta^*_n,x)|
\leq
\limsup_{n\to \infty}
\sup_{x\in\mathbb{R}}
|\Phi_n(\beta^*_n,x)|<\infty,
\]
with probability one.
\item
For $i=1,2$,
\[
\sup_{x\in\mathbb{R}}
\sup_{|\beta-\beta_0| \leq \varepsilon} |D^{(i)}(\beta,x)|<\infty.
\]
\item
For $i=1,2$ and for any sequence $\beta^*_n$, such that $\beta^*_n\to \beta_0$ almost surely,
\[
\limsup_{n\to \infty}
\sup_{x\in\mathbb{R}}
|D_n^{(i)}(\beta^*_n,x)|<\infty,
\]
with probability one.
\end{enumerate}
\end{lemma}
\noindent
The proof can be found in the appendix.

Obviously, we will approximate $\Phi_n(\hb,x)$ and $\Phi_n(\beta_0,x)$
by $\Phi(\beta_0,x)$.
According to the law of large numbers, $\Phi_n$ will converge to $\Phi$, for $\beta$ and $x$ fixed.
However, we need uniform convergence at proper rates.
This is established by the following lemma.

\begin{lemma}
\label{lemma:conv Phin}
Suppose that condition (A2) holds and $\hb \to \beta_0$, with probability one.
Then,
\[
\sup_{x \in \mathbb{R}}
\left|\Phi_n(\hb,x)-\Phi(\beta_0,x)\right|\to 0,
\]
with probability one.
Moreover,
\begin{equation}
\label{eq:rate Phin}
\sqrt{n}\sup_{x\in\mathbbm{R}}\left| \Phi_n(\beta_0,x)-\Phi(\beta_0,x) \right|=\O_p(1).
\end{equation}
\end{lemma}
\begin{proof}
For all $x\in \mathbb{R}$, write
\[
|\Phi_n(\hb,x)-\Phi(\beta_0,x)|
\leq
|\Phi_n(\hb,x)-\Phi_n(\beta_0,x)|
+
|\Phi_n(\beta_0,x)-\Phi(\beta_0,x)|.
\]
For the second term on the right hand side, consider the class of functions
\[
\mathcal{G}=\left\{g(u,z;x): x\in\mathbb{R} \right\},
\]
where for each $x\in\mathbb{R}$ and $\beta_0\in\mathbbm{R}^p$ fixed,
\[
g(u,z;x)=\{u\geq x\}\exp(\beta_0'z)
\]
is a product of an indicator and a fixed function.
It follows that $\mathcal{G}$ is a VC-subgraph class (e.g., see Lemma 2.6.18 in~\cite{vaartwellner})
and its envelope $G=\exp(\beta_0'z)$ is square integrable under condition~(A2).
Standard results from empirical process theory~\cite{vaartwellner} yield that the class of functions $\mathcal{G}$ is Glivenko-Cantelli,
i.e.,
\begin{equation}
\label{eq:tsiatis}
\sup_{x \in \mathbb{R}}
|\Phi_n(\beta_0,x)-\Phi(\beta_0,x)|
=
\sup_{g \in \mathcal{G}}\left| \int g(u,z;x)\,\mathrm{d}(\mathbbm{P}_n-P)(u,\delta,z) \right|\to 0,
\end{equation}
with probability one.
Moreover, $\mathcal{G}$ is a Donsker class, i.e.,
\[
\sqrt{n}\int g(u,z;x)\,\mathrm{d}(\mathbbm{P}_n-P)(u,\delta,z)=\O_p(1),
\]
so that~\eqref{eq:rate Phin} follows by continuous mapping theorem.
Finally, by Taylor expansion and the Cauchy-Schwarz inequality, it follows that
\[
\sup_{x \in \mathbb{R}}
\left|\Phi_n(\hb,x)-\Phi_n(\beta_0,x)\right|
=
\sup_{x \in \mathbb{R}}
\left| (\hb-\beta_0)'D_n^{(1)}(\beta^*,x)\right|
\leq
|\hb-\beta_0|
\sup_{x \in \mathbb{R}}
\left| D_n^{(1)}(\beta^*,x)\right|,
\]
for some $\beta^*$, for which $|\beta^*-\beta_0|\leq |\hb-\beta_0|$.
Together with~\eqref{eq:tsiatis},
from the strong consistency of $\hb$ (e.g., see Theorem~3.1 in~\cite{tsiatis:1981}) and Lemma~\ref{lemma:bounds Phin},
the lemma follows.
\end{proof}
The previous results can be used to prove a first step in the direction of proving Theorem~\ref{theorem:consistency_estimators},
i.e., suitable uniform approximation of $\Lambda_n$ and $F_n$ by $\Lambda_0$ and $F_0$.
Strong uniform consistency of $\Lambda_n$ and process convergence of $\sqrt{n}(\Lambda_n-\Lambda_0)$
has been established in~\cite{kosorok}, under the stronger assumption of bounded covariates.
Weak consistency has been derived or mentioned before, see for example~\cite{prenticekalbfleisch:2003}.
\begin{theorem}
\label{theorem:conv_breslow}
Under the assumptions (A1) and (A2), for all $0<M<\tau_H$,
\[
\sup_{x\in[0,M]}\left|\BR(x)-\Lambda_0(x)\right|\to 0,
\]
with probability one and
$\sqrt{n} \sup_{x\in[0,M]} \left| \BR(x)-\Lambda_0(x) \right|=\O_p(1)$.
\end{theorem}
\begin{proof}
From the expression for the baseline cumulative hazard function in~\eqref{eq:relation cum base haz}
together with~\eqref{def:Phi} and~\eqref{eq:relation Breslow Phi_n}, it follows that
\[
\begin{split}
\sup_{x\in[0,M]}\left|\BR(x)-\Lambda_0(x)\right|
&\leq
\sup_{x\in[0,M]} \left| \int
\delta\{u\leq x\}
\left( \frac{1}{\Phi_n(\hb,u)}-\frac{1}{\Phi_n(\beta_0,u)} \right) \ \mathrm{d}\mathbbm{P}_n(u,\delta,z)\right|\\
&\quad+
\sup_{x\in[0,M]}\left| \int
\delta\{u\leq x\}
\left( \frac{1}{\Phi_n(\beta_0,u)}-\frac{1}{\Phi(\beta_0,u)} \right)\ \mathrm{d}\mathbbm{P}_n(u,\delta,z)\right|\\
&\quad+
\sup_{x\in[0,M]}\left| \int
\frac{\delta\{u\leq x\}}{\Phi(\beta_0,u)}\  \mathrm{d}\left(\mathbbm{P}_n-P\right)(u,\delta,z)\right|\\
&=A_n+B_n+C_n.
\end{split}
\]
Starting with the first term on the right hand side, note that
\begin{equation}
\label{eq:bound A_n}
A_n
\leq
\frac{|\hb-\beta_0|}{\Phi_n(\hb,M)\Phi_n(\beta_0,M)}
\sup_{x\in\mathbbm{R}}\left| D_n^{(1)}(\beta^*,x)\right|
\end{equation}
for some $|\beta^*-\beta_0|\leq |\hb-\beta_0|$.
According to Lemma~\ref{lemma:bounds Phin}, the right hand side is bounded by $C|\hb-\beta_0|$,
for some $C>0$.
Since $\hb$ is strong consistent and $|\hb-\beta_0|=\O_p(n^{-1/2})$,
(e.g., see Theorems 3.1 and 3.2 in~\cite{tsiatis:1981}),
it follows that $A_n\to0$ almost surely and
$A_n=\O_p(n^{-1/2})$.
Similarly,
\begin{equation}
\label{eq:bound Bn}
B_n\leq \frac{1}{\Phi_n(\beta_0,M)\Phi(\beta_0,M)}\sup_{x\in\mathbbm{R}}\left| \Phi_n(\beta_0,x)-\Phi(\beta_0,x) \right|.
\end{equation}
From Lemmas~\ref{lemma:bounds Phin} and~\ref{lemma:conv Phin},
it follows that $B_n\to0$ almost surely and $B_n=\O_p(n^{-1/2})$.
For the last term $C_n$, consider the class of functions
$\mathcal{H}=\left\{ h(u,\delta;x):x\in[0,M] \right\}$,
where for each $x\in[0,M]$, with $M<\tau_H$ and $\beta_0\in\mathbbm{R}^p$ fixed,
\[
h(u,\delta;x)=\frac{\delta\{u\leq x\}}{\Phi(\beta_0,u)}.
\]
The function $h$ is a product of indicators and a fixed uniformly bounded monotone function.
Similar to the arguments given in the proof of Lemma~\ref{lemma:conv Phin}, it follows that the class~$\mathcal{H}$ is Glivenko-Cantelli,
i.e.,
\[
\sup_{h\in\mathcal{H}}\left| \int h(u,\delta;\cdot)\mathrm{d}(\mathbbm{P}_n-P)(u,\delta,z) \right|\to 0,
\]
almost surely, which gives the first statement of the lemma.
Moreover, $\mathcal{H}$ is a Donsker class and hence the second statement of the lemma follows by continuous mapping theorem. This completes the proof.
\end{proof}
Strong uniform consistency of $F_n$ follows immediately from the strong consistency of the Breslow estimator
established in Theorem~\ref{theorem:conv_breslow}, and is stated in the next corollary.
\begin{corollary}
\label{cor:consistency_F_n}
Under the assumptions (A1) and (A2) and for all $0<M<\tau_H$,
\[
\sup_{x\in[0,M]}|F_n(x)-F_0(x)|\to 0,
\]
with probability one.
\end{corollary}
\begin{proof}
The proof is straightforward and follows immediately from Theorem~\ref{theorem:conv_breslow},
relations~\eqref{eq:relation F0} and~\eqref{eq: relation Fn}, together
with the fact that $|\text{e}^{-y}-1|\leq 2|y|$, as $y\to0$.
\end{proof}
Note that the estimators in Theorem~\ref{theorem:consistency_estimators} of the baseline hazard are essentially
the slopes of the GCM of $V_n$.
For this reason, as a final preparation for the proof of Theorem~\ref{theorem:consistency_estimators},
we establish uniform convergence of the GCM of $V_n$ by the following lemma.
This lemma is completely similar to Lemma~4.3 in~\cite{huangwellner:1995}.
Its proof can be found in the appendix.
\begin{lemma}
\label{lemma:conv Wn Vn hat}
Assume that $\Lambda_0$ is convex on $[0,\tau_H]$ and that conditions (A1) and (A2) hold.
Let $\hb$ be the maximum partial likelihood estimator and define
\begin{equation}
\label{def:Wn_hat}
\widehat{W}_n(x)=W_n(\hb,x)-W_n(\hb,T_{(1)}),
\qquad
x\geq T_{(1)},
\end{equation}
where $W_n$ is defined in~\eqref{def:W_n}.
Let $\big(\widehat{W}_n(x),\widehat{V}_n(x)\big)$ be the GCM of
$\big( \widehat{W}_n(x),V_n(x) \big)$, for~$x\in[T_{(1)},T_{(n)}]$,
where $V_n$ is defined in~\eqref{def:V_n}.
Then
\begin{equation}
\label{eq:conv sup V_n hat}
\sup_{x\in [T_{(1)},T_{(n)}]} \left| \widehat{V}_n(x)-V(x) \right| \to 0,
\end{equation}
with probability one, where $V(x)=H^{uc}(x)$,
as defined just below~\eqref{eq:relation base haz}.
\end{lemma}

Obviously, in the nonincreasing case, similar to~\eqref{eq:conv sup V_n hat} one can show
\begin{equation}
\label{eq:conv Y_n}
\sup_{x\in [0,T_{(n)}]}
\left|\widehat{Y}_n(x)-V(x) \right|\to 0,
\end{equation}
almost surely, where $\big(W_n(\hb,x),\widehat{Y}_n(x)\big)$
is the LCM of $\big(W_n(\hb,x),Y_n(x) \big)$,
with $Y_n$ defined in~\eqref{def:Y_n}.
We are now in the position to prove Theorem~\ref{theorem:consistency_estimators},
which establishes strong consistency of the estimators.
\begin{proof}[Proof of Theorem~\ref{theorem:consistency_estimators}]
First consider the second statement of case~(i).
Since $\GCM$ is convex on the open interval $(0,\tau_H)$,
it admits in every point $x_0\in(0,\tau_H)$ a finite left and a right derivative, denoted by $\widetilde{\Lambda}_{n}^{-}$ and $\widetilde{\Lambda}_{n}^{+}$ respectively. Moreover, for any fixed $x_0\in(0,\tau_H)$ and for sufficiently small $\delta>0$, it follows that
\[
\frac{\GCM(x_0)-\GCM(x_0-\delta)}{\delta} \leq \widetilde{\Lambda}_{n}^{-}(x_0)
\leq \widetilde{\Lambda}_{n}^{+}(x_0) \leq \frac{\GCM(x_0+\delta)-\GCM(x_0)}{\delta}.
\]
When $n\to\infty$, then for any $0<M<\tau_H$,
\begin{equation}
\label{eq:marshall's lemma}
\sup_{x\in[0,M]}\left| \GCM(x)-\Lambda_0(x) \right| \leq \sup_{x\in[0,M]}\left|\BR(x)-\Lambda_0(x) \right|.
\end{equation}
This is a variation of Marshall's lemma and can be proven similar to (7.2.3) in~\cite{robertsonetal}
or Lemma~4.1 in~\cite{huangwellner:1995}.
By convexity of $\Lambda_0$ and the fact that $\GCM$ is the greatest convex function below $\BR$,
one must have
\[
\Lambda_0(x)-\delta_n \leq \GCM(x) \leq \BR(x),
\]
where $\delta_n=\sup_{x\in[0,M]} |\Lambda_0(x)-\BR(x)|$,
which yields inequality~\eqref{eq:marshall's lemma}.
From~\eqref{eq:marshall's lemma} and Theorem~\ref{theorem:conv_breslow},
by first letting $n\rightarrow \infty$ and then $\delta\rightarrow0$, we find
\[
\lambda_0(x_0-)\leq \liminf_{n\rightarrow \infty} \widetilde{\Lambda}_{n}^{-}(x_0) \leq \limsup_{n\rightarrow \infty} \widetilde{\Lambda}_{n}^{-}(x_0) \leq \limsup_{n\rightarrow \infty} \widetilde{\Lambda}_{n}^{+}(x_0) \leq \lambda_0(x_0+).
\]
Because $\ho(x_0)=\widetilde{\Lambda}_{n}^{-}(x_0)$, this proves that $\ho$ is a strong consistent estimator.

For $\hoMLE$, first note that since $\widehat{V}_n$
is convex on the open interval $(0,\tau_H)$, it admits in every point $x_0\in(0,\tau_H)$ a finite left and a right derivative, denoted by $\widehat{V}_{n}^{-}$ and $\widehat{V}_{n}^{+}$ respectively, where
\[\begin{split}
\widehat{V}_{n}^{-}(x)
&=
\lim_{\delta\downarrow0}\frac{\widehat{V}_n(x)-\widehat{V}_n(x-\delta)} {\widehat{W}_n(x)-\widehat{W}_n(x-\delta)},\\
\widehat{V}_{n}^{+}(x)
&=
\lim_{\delta\downarrow0}\frac{\widehat{V}_n(x+\delta)-\widehat{V}_n(x)} {\widehat{W}_n(x+\delta)-\widehat{W}_n(x)}.
\end{split}\]
For any fixed $x\in(0,\tau_H)$ and for sufficiently small $\delta>0$, it follows that
\[
\frac{\widehat{V}_n(x_0)-\widehat{V}_n(x_0-\delta)}{\widehat{W}_n(x_0)-\widehat{W}_n(x_0-\delta)}
\leq
\widehat{V}_{n}^{-}(x_0)
\leq
\widehat{V}_{n}^{+}(x_0)
\leq
\frac{\widehat{V}_n(x_0+\delta)-\widehat{V}_n(x_0)}{\widehat{W}_n(x_0+\delta)-\widehat{W}_n(x_0)}.
\]
If we define
\[
W_0(x)
=
\int_0^x \Phi(\beta_0,s)\,\text{d}s,
\]
then by making use of Lemma~\ref{lemma:conv Wn Vn hat}, together with
\begin{equation}
\label{eq:bound Wnhat-W0 main text}
\sup_{x\in[T_{(1)},T_{(n)}]} \left| \widehat{W}_n(x)-W_0(x) \right|
\leq \tau_H \sup_{x \in\mathbbm{R}} \left| \Phi_n(\hb,x)-\Phi(\beta_0,x) \right| \to 0,
\end{equation}
with probability one (see the proof of Lemma~\ref{lemma:conv Wn Vn hat} in the appendix) and letting $n\to \infty$,
we obtain
\[
\frac{V(x_0)-V(x_0-\delta)}{W_0(x_0)-W_0(x_0-\delta)}
\leq
\liminf_{n\rightarrow \infty} \widehat{V}_n^-(x_0)
\leq
\limsup_{n\rightarrow \infty} \widehat{V}_n^+(x_0)
\leq
\frac{V(x_0+\delta)-V(x_0)}{W_0(x_0+\delta)-W_0(x_0)}.
\]
Furthermore, by letting $\delta\to 0$, together with
the fact that, according to~\eqref{eq:relation base haz} and~\eqref{def:Phi},
$\lambda_0$ can also be represented as
\[
\lambda_0(x)=\frac{\mathrm{d} V(x)/\mathrm{d}x}{\mathrm{d}W_0(x)/\mathrm{d}x},
\]
we get
\[
\lambda_0(x_0-)
\leq
\liminf_{n\rightarrow \infty} \widehat{V}_n^-(x_0)
\leq
\limsup_{n\rightarrow \infty} \widehat{V}_n^-(x_0)
\leq
\limsup_{n\rightarrow \infty} \widehat{V}_n^+(x_0)
\leq
\lambda_0(x_0+),
\]
which completes the proof of (i), since $\hoMLE(x_0)=\widehat{V}_n^-(x_0)$.
The proofs of (ii) and (iii) are completely analogous,
using~\eqref{eq:conv Y_n} and Corollary~\ref{cor:consistency_F_n}.
\end{proof}

\section{Inverse processes}
\label{sec:limit_distr_preliminary}
\noindent
To obtain the limit distribution of the estimators, we follow the approach proposed in~\cite{groeneboom:1985}.
For each proposed estimator, we define an inverse process and establish its asymptotic distribution.
The asymptotic distribution of the estimators then emerges via the switching relationships.
The inverse processes are defined in terms of some local processes and this section
is devoted to acquire the weak convergence of these local processes.
Furthermore, the inverse processes need to be bounded in probability.
This result, along with the limiting distribution of the inverse processes and hence of the estimators are deferred to
Section~\ref{sec:limit distribution}.

In order to keep the exposition brief, we do not treat all five separate cases in detail,
but we confine ourselves to the most important ones, as the other cases can be handled similarly.
In the case of a nondecreasing $\lambda_0$, the distribution of the NPMLE $\hoMLE$
can be obtained through the study of the inverse process
\begin{equation}
\label{def:inverse MLE}
\widehat{U}_n^{\lambda}(a)=\argmin_{x\in[T_{(1)},T_{(n)}]}
\left\{ V_n(x)-a\widehat{W}_n(x)\right\},
\end{equation}
for $a>0$,
where  $V_n$ and $\widehat{W}_n$ have been defined in~\eqref{def:V_n} and~\eqref{def:Wn_hat}.
Succeedingly,
for a given $a>0$, the switching relationship holds,
i.e., $\widehat{U}_n^{\lambda}(a)\geq x$ if and only if $\hoMLE(x)\leq a$
with probability one, so that after scaling, it follows that
\begin{equation}
\label{eq:switching incr}
n^{1/3}
\left[
\hoMLE(x_0)-\lambda_0(x_0)
\right]>a
\Leftrightarrow
n^{1/3}
\left[
\widehat{U}_n^{\lambda}(\lambda_0(x_0)+n^{-1/3}a)-x_0
\right]<0,
\end{equation}
for $0<x_0<\tau_H$, with probability one.
A similar relationship holds for $\ho$ and the corresponding inverse process
\begin{equation}
\label{def:inverse LS haz}
\widetilde{U}_n^{\lambda}(a)=
\argmin_{x\in[0,T_{(n)}]}
\left\{\BR(x)-ax\right\}.
\end{equation}
For the nonincreasing density estimator $\fo$,
we consider the inverse process
\begin{equation}
\label{def:inverse_process_f}
\widetilde{U}_n^f(a)=\argmax_{x\in[0,T_{(n)}]}\left\{ F_n(x)-ax \right\},
\end{equation}
where argmax denotes the largest location of the maximum.
In this case, instead of~\eqref{eq:switching incr}, we have
\begin{equation}
\label{eq:switching decr}
n^{1/3}
\left[
\fo(x_0)-f_0(x_0)
\right]>a
\Leftrightarrow
n^{1/3}
\left[
\widetilde{U}_n^{f}(f_0(x_0)+n^{-1/3}a)-x_0
\right]>0,
\end{equation}
Similarly, in the case of estimating a nonincreasing $\lambda_0$,
we consider inverse processes $\widehat{U}_n^{\lambda}$ and~$\widetilde{U}_n^{\lambda}$
defined with argmax instead of argmin in~\eqref{def:inverse MLE}
and~\eqref{def:inverse LS haz} and we have switching relations similar to~\eqref{eq:switching decr}.

From the definition of the inverse process in~\eqref{def:inverse LS haz} and given that the argmin is invariant
under addition of and multiplication with positive constants, it can be derived that
\begin{equation}
\label{def:Un scaled haz LS}
n^{1/3}
\left[
\widetilde{U}_n^{\lambda}(\lambda_0(x_0)+n^{-1/3}a)-x_0
\right]
=
\argmin_{x\in I_n(x_0)}
\left\{
\widetilde{\mathbbm{Z}}_n^{\lambda}(x)-ax
\right\}
\end{equation}
where $I_n(x_0)=[-n^{1/3}x_0,n^{1/3}(T_{(n)}-x_0)]$ and
\begin{equation}
\label{def:Zn incr haz LS}
\widetilde{\mathbbm{Z}}_n^{\lambda}(x)
=
n^{2/3}
\left[\BR(x_0+n^{-1/3}x)-\BR(x_0)-n^{-1/3}\lambda_0(x_0) x\right].
\end{equation}
Likewise, $n^{1/3}\left[ \widehat{U}_n^{\lambda}(\lambda_0(x_0)+n^{-1/3}a)-x_0\right]$ is equal to
\begin{equation}
\label{def:Un scaled haz ML}
\argmin_{x\in I_n'(x_0)}\left\{ \widehat{\mathbbm{Z}}_n^{\lambda}(x)- \frac{n^{1/3}a}{\Phi(\beta_0,x_0)} \left[ \widehat{W}_n(x_0+n^{-1/3}x)-\widehat{W}_n(x_0) \right]\right\},
\end{equation}
where $I_n'(x_0)=[-n^{1/3}(x_0-T_{(1)}),n^{1/3}(T_{(n)}-x_0)]$ and
\begin{equation}
\label{def:Zn incr haz MLE}
\begin{split}
\widehat{\mathbbm{Z}}_n^{\lambda}(x)
=
\frac{n^{2/3}}{\Phi(\beta_0,x_0)}
\Bigg(
V_n(x_0&+n^{-1/3}x)-V_n(x_0)\\
&-\lambda_0(x_0)\left[\widehat{W}_n(x_0+n^{-1/3}x)-\widehat{W}_n(x_0)\right]
\Bigg),
\end{split}
\end{equation}
and similarly
\begin{equation}
\label{def:Un dens}
n^{1/3}
\left[\widetilde{U}_n^f(f_0(x_0)+n^{-1/3}a)-x_0\right]
=
\argmax_{x\in I_n(x_0)}\{\mathbbm{\widetilde{Z}}^f_n(x)-ax\},
\end{equation}
where
\begin{equation}
\label{def:Zn dens}
\mathbbm{\widetilde{Z}}^f_n(x)
=
n^{2/3}
\left[
F_n(x_0+n^{-1/3}x)-F_n(x_0)-n^{-1/3}f_0(x_0) x
\right].
\end{equation}
In the case of estimating a nonincreasing $\lambda_0$, we consider
the argmax of the processes~\eqref{def:Zn incr haz MLE} and~\eqref{def:Zn incr haz LS}.
Before investigating the asymptotic behavior of the above processes, we
first need to establish the following technical lemma.
It provides a sufficient bound on the order of shrinking increments of
an empirical process that we will encounter later on.
\begin{lemma}
\label{lemma:Jon}
Assume (A1) and (A2).
Let $x_0\in(0,\tau_H)$ fixed and suppose that
\begin{equation}
\label{eq:cond H^uc}
H^{uc}
\text{ is continuously differentiable in a neighborhood of $x_0$.}
\end{equation}
Then, for any $k=1,2,\ldots$,
\[
\sup_{|x|\leq k}
\left|
\int
\delta
\left(
\{u\leq x_0+n^{-1/3}x\}-\{u\leq x_0\}
\right)
\left( \frac{1}{\Phi_n(\beta_0,u)}-\frac{1}{\Phi(\beta_0,u)} \right)\mathrm{d}(\mathbbm{P}_n-P)(u,\delta,z)\right|
\]
is of the order $\O_p(n^{-7/6}\log n)$.
\end{lemma}
\begin{proof}
Take $0\leq x\leq k$ and consider the class of functions
\begin{equation}
\label{def:class_F_n}
\mathcal{F}_n=\left\{ f_n(u,\delta,z;x): 0\leq x\leq k \right\},
\end{equation}
where for each $0\leq x\leq k$,
\[
f_n(u,\delta,z;x)
=
\delta
\{x_0< u\leq x_0+n^{-1/3}x\}
\left( \frac{1}{\Phi_n(\beta_0,u)}-\frac{1}{\Phi(\beta_0,u)}\right).
\]
Correspondingly, consider the class $\mathcal{G}_{n,k,\alpha}$ consisting of functions
\[
g(u,\delta,z;y,\Psi)=
\delta\{x_0< u\leq x_0+y\}
\left( \frac{1}{\Psi(u)}-\frac{1}{\Phi(\beta_0,u)} \right).
\]
where $0\leq y \leq n^{-1/3}k$ and $\Psi$ is nonincreasing left continuous, such that
\[
\Psi(x_0+n^{-1/3}k)
\geq
K
\quad\text{and}\quad
\sup_{u\in \mathbbm{R}}\left| \Psi(u)-\Phi(\beta_0,u) \right|\leq \alpha,
\]
where $K=\Phi(\beta_0,(x_0+\tau_H)/2)/2$.
Then, for any $\alpha>0$ and $k=1,2,\ldots$,
\[
P\left( \mathcal{F}_n\subset \mathcal{G}_{n,k,\alpha} \right) \to 1,
\]
by Lemma~\ref{lemma:conv Phin}.
Furthermore, the class $\mathcal{G}_{n,k,\alpha}$ has envelope
\[
G(u,\delta,z)=
\delta\{x_0< u\leq x_0+n^{-1/3}k\}
\frac{\alpha}{K^2},
\]
for which it follows from~\eqref{eq:cond H^uc}, that
\[
\| G \|_{P,2}^2
=
\int G(u,\delta,z)^2\,\mathrm{d}P(u,\delta,z)
=
\frac{\alpha^2}{K^4}P(x_0< T\leq x_0+n^{-1/3}k,\Delta=1)
=
\O(\alpha^2 kn^{-1/3}).
\]
Since the functions in $\mathcal{G}_{n,k,\alpha}$ are sums and products of bounded monotone functions,
its entropy with bracketing satisfies
\[
\log N_{[\,]}(\varepsilon,\mathcal{G}_{n,k,\alpha},L_2(P))\lesssim \frac{1}{\varepsilon},
\]
see e.g., Theorem~2.7.5 in~\cite{vaartwellner} and Lemma~9.25 in~\cite{kosorok},
and hence, for any $\delta>0$, the bracketing integral
\[
J_{[\,]}(\delta,\mathcal{G}_{n,k,\alpha},L_2(P))=\int_0^\delta\sqrt{1+\log N_{[\,]}(\varepsilon\|G\|_2,\mathcal{G}_{n,k,\alpha},L_2(P))}\,\mathrm{d}\varepsilon
<\infty.
\]
By Theorem~2.14.2 in~\cite{vaartwellner}, we have
\[
\begin{split}
\mathbb{E}\left\|\sqrt{n}\int g(u,\delta,z;y,\Psi)\mathrm{d}(\mathbbm{P}_n-P)(u,\delta,z) \right\|_{\mathcal{G}_{n,k,\alpha}}
&\leq
J_{[\,]}(1,\mathcal{G}_{n,k,\alpha},L_2(P))\|G\|_{P,2}\\
&=
\O(\alpha k^{1/2}n^{-1/6}),
\end{split}
\]
where $\|\cdot\|_{\mathcal{F}}$ denotes the supremum over the class of functions $\mathcal{F}$.
Now, according to~\eqref{eq:rate Phin}
\[
(\log n)^{-1}\sqrt{n}\sup_{x\in\mathbbm{R}}\left| \Phi_n(\beta_0,x)-\Phi(\beta_0,x) \right|\to 0,
\]
in probability.
Therefore, if we choose $\alpha=n^{-1/2}\log n$, this gives
\[
\mathbb{E}
\left\|
\int g(u,\delta,z;y,\Psi)\mathrm{d}(\mathbbm{P}_n-P)(u,\delta,z)
\right\|_{\mathcal{G}_{n,k,\alpha}}
=
\O(k^{1/2} n^{-7/6}\log n)
\]
and hence by the Markov inequality,
this proves the lemma for the case $0\leq x\leq k$.
The argument for $-k\leq x\leq 0$ is completely similar.
\end{proof}
Our approach in deriving the asymptotic distribution of the monotone estimators
involves application of results from~\cite{kimpollard:1990}.
To this end, we first determine the limiting processes of~\eqref{def:Zn incr haz MLE}, \eqref{def:Zn incr haz LS}
and~\eqref{def:Zn dens}.
\begin{lemma}
\label{lemma:relation_Z_n_tilde_hat}
Suppose that (A1) and (A2) hold. Assume~\eqref{eq:cond H^uc} and that
\begin{equation}
\label{eq:cond lambda0}
\text{$\lambda_0$ is continuously differentiable in a neighborhood of $x_0$}.
\end{equation}
Moreover, assume that
\begin{equation}
\label{eq:cond Phi}
x\mapsto \Phi(\beta_0,x)
\text{ is continuously differentiable in a neighborhood of $x_0$}.
\end{equation}
Then, for any $k=1,2,\ldots$,
\[
\sup_{|x|\leq k} \left| \widetilde{\mathbbm{Z}}_n^{\lambda}(x)- \widehat{\mathbbm{Z}}_n^{\lambda}(x) \right|\to 0,
\]
in probability, where the processes $\widetilde{\mathbbm{Z}}_n^{\lambda}$ and $\widehat{\mathbbm{Z}}_n^{\lambda}$
are defined in~\eqref{def:Zn incr haz LS} and~\eqref{def:Zn incr haz MLE}, respectively.
\end{lemma}
\begin{proof}
We will prove that for any $k=1,2,\ldots,$
\[
\sup_{x\in[0,k]}\left| \widetilde{\mathbbm{Z}}_n^{\lambda}(x)-\widehat{\mathbbm{Z}}_n^{\lambda}(x)  \right|\to 0,
\]
in probability, since the result for $-k\leq x\leq 0$ follows completely analogous.
Write
\[\begin{split}
&
\Phi(\beta_0,x_0)\left( \widetilde{\mathbbm{Z}}_n^{\lambda}(x)-\widehat{\mathbbm{Z}}_n^{\lambda}(x) \right)\\
&=
n^{2/3} \int
\delta\left\{x_0<u\leq x_0+n^{-1/3}x\right\}
\left( \frac{\Phi(\beta_0,x_0)}{\Phi_n(\hb,u)}-1 \right)\mathrm{d}\mathbbm{P}_n(u,\delta,z)\\
&\qquad-
n^{2/3}\lambda_0(x_0)\int_{x_0}^{x_0+n^{-1/3}x}
\left[ \Phi(\beta_0,x_0)-\Phi_n(\hb,s) \right]\mathrm{d}s\\
&=
n^{2/3} \int
\delta\left\{x_0<u\leq x_0+n^{-1/3}x\right\}
\left( \frac{\Phi(\beta_0,x_0)}{\Phi_n(\hb,u)}-\frac{\Phi(\beta_0,x_0)}{\Phi_n(\beta_0,u)} \right)\mathrm{d}\mathbbm{P}_n(u,\delta,z)\\
&\qquad+
n^{2/3} \int
\delta\left\{x_0<u\leq x_0+n^{-1/3}x\right\}
\left( \frac{\Phi(\beta_0,x_0)}{\Phi_n(\beta_0,u)}-1 \right)\mathrm{d}\mathbbm{P}_n(u,\delta,z)\\
&\qquad-
n^{2/3}\lambda_0(x_0)\int_{x_0}^{x_0+n^{-1/3}x}
\left[ \Phi(\beta_0,x_0)-\Phi_n(\beta_0,s) \right]\mathrm{d}s\\
&\qquad-
n^{2/3}\lambda_0(x_0)\int_{x_0}^{x_0+n^{-1/3}x}
\left[ \Phi_n(\beta_0,s)-\Phi_n(\hb,s) \right]\mathrm{d}s\\
&=
A_{n1}(x)+A_{n2}(x)+A_{n3}(x)+A_{n4}(x).
\end{split}\]
We will show that the supremum of all four terms on the right hand side tend to zero in probability.
Similar to~\eqref{eq:bound A_n}, according to Lemma~\ref{lemma:bounds Phin},
\[
|A_{n1}(x)|
\leq
C|\hb-\beta_0|
n^{2/3}\int
\left\{x_0<u\leq x_0+n^{-1/3}x\right\}\,\mathrm{d}\mathbbm{P}_n(u,\delta,z),
\]
for some $C>0$.
Since, $|\hb-\beta_0|=\O_p(n^{-1/2})$ and
\[
\int\left\{x_0<u\leq x_0+n^{-1/3}x\right\}\,\mathrm{d}(\mathbbm{P}_n-P)(u,\delta,z)
=
\O_p(n^{-2/3}x^{1/2})+\O_p(n^{-1/3}x),
\]
it follows that
\begin{equation}
\label{eq:bound An1}
|A_{n1}(x)|=\O_p(n^{-1/2}x^{1/2})+\O_p(n^{-1/6}x),
\end{equation}
and likewise, $|A_{n4}(x)|=\O_p(n^{-1/6}x)$.
Furthermore, write
\[\begin{split}
A_{n2}(x)
&=
n^{2/3}
\int\delta\left\{x_0<u\leq x_0+n^{-1/3}x\right\}
\left( \frac{\Phi(\beta_0,x_0)}{\Phi_n(\beta_0,u)}-\frac{\Phi(\beta_0,x_0)}{\Phi(\beta_0,u)} \right)\mathrm{d}\left( \mathbbm{P}_n-P \right)(u,\delta,z)\\
&\qquad+
n^{2/3}
\int\delta\left\{x_0<u\leq x_0+n^{-1/3}x\right\}
\left(\frac{\Phi(\beta_0,x_0)}{\Phi(\beta_0,u)}-1\right)\mathrm{d}\left(\mathbbm{P}_n -P \right)(u,\delta,z)\\
&\qquad+
n^{2/3}
\int\delta\left\{x_0<u\leq x_0+n^{-1/3}x\right\}
\left( \frac{\Phi(\beta_0,x_0)}{\Phi_n(\beta_0,u)} -\frac{\Phi(\beta_0,x_0)}{\Phi(\beta_0,u)}\right) \mathrm{d}P(u,\delta,z)\\
&\qquad+
n^{2/3}
\int\delta\left\{x_0<u\leq x_0+n^{-1/3}x\right\}
\left( \frac{\Phi(\beta_0,x_0)}{\Phi(\beta_0,u)}-1 \right)\mathrm{d}P(u,\delta,z)\\
&=
B_{n1}(x)+B_{n2}(x)+B_{n3}(x)+B_{n4}(x).
\end{split}\]
According to Lemma~\ref{lemma:Jon},
\begin{equation}
\label{eq:bound Bn1}
\sup_{0\leq x\leq k}
\left| B_{n1}(x) \right|=\O_p(n^{-1/2}\log n).
\end{equation}
For the term $B_{n2}$, consider the class $\mathcal{F}$ consisting of functions
\[
f(u,\delta,z;x)
=
\delta\{x_0< u\leq x_0+n^{-1/3}x\}\left( \frac{\Phi(\beta_0,x_0)}{\Phi(\beta_0,u)}-1 \right),
\]
where $0\leq x\leq k$, with envelope
\[
F(u)=\delta\{x_0< u\leq x_0+n^{-1/3}k\}
\left( \frac{\Phi(\beta_0,x_0)}{\Phi(\beta_0,x_0+n^{-1/3}k)}-1 \right).
\]
Then, the $L_2(P)$ norm of the envelope satisfies
\[
\|F\|_{P,2}^2
=
\left( \frac{\Phi(\beta_0,x_0)}{\Phi(\beta_0,x_0+n^{-1/3}k)}-1 \right)^2
\left[ H^{uc}(x_0+n^{-1/3}k)-H^{uc}(x_0) \right]
=
\O(n^{-1}),
\]
according to~\eqref{eq:cond H^uc} and Lemma~\ref{lemma:bounds Phin},
so that by arguments similar as in the proof of Lemma~\ref{lemma:Jon},
\begin{equation}
\label{eq:bound Bn2}
\sup_{0\leq x\leq k}|B_{n2}(x)|=\O_p(n^{-1/3}).
\end{equation}
For the term $B_{n3}$,
similar to the treatment of the right hand side of~\eqref{eq:bound Bn},
it follows that
\begin{equation}
\label{eq:bound Bn3}
|B_{n3}(x)|
\leq
n^{2/3}\O_p(n^{-1/2})
\left| H^{uc}(x_0+n^{-1/3}x)-H^{uc}(x_0) \right|
=
\O_p(n^{-1/6}x),
\end{equation}
by condition~\eqref{eq:cond H^uc}.
Next, we combine $B_{n4}(x)$ with $A_{n3}(x)$.
First write
\[
\begin{split}
A_{n3}(x)
&=
n^{2/3}\lambda_0(x_0)\int_{x_0}^{x_0+n^{-1/3}x}
\left[\Phi_n(\beta_0,s)-\Phi(\beta_0,s) \right]\,\mathrm{d}s\\
&\qquad+
n^{2/3}\lambda_0(x_0)\int_{x_0}^{x_0+n^{-1/3}x}
[\Phi(\beta_0,s)-\Phi(\beta_0,x_0) ]\,\mathrm{d}s\\
&=
C_{n1}(x)+C_{n2}(x).
\end{split}
\]
As for $C_{n1}$,
\begin{equation}
\label{eq:bound Cn1}
|C_{n1}(x)|
\leq
n^{1/3}x \lambda_0(x_0) \sup_{x\in\mathbbm{R}}\left| \Phi_n(\beta_0,x)-\Phi(\beta_0,x) \right|
=\O_p(n^{-1/6}x),
\end{equation}
according to Lemma~\ref{lemma:conv Phin}.
Finally, using~\eqref{eq:relation base haz} and~\eqref{def:Phi},
\begin{equation}
\label{eq:bound Bn4+Cn2}
\begin{split}
B_{n4}(x)+C_{n2}(x)
&=
n^{2/3}
\int_{x_0}^{x_0+n^{-1/3}x}
\left[
\Phi(\beta_0,x_0)-\Phi(\beta_0,u)
\right]
\lambda_0(u)\,\mathrm{d}u \\
&\qquad+
n^{2/3}\lambda_0(x_0)\int_{x_0}^{x_0+n^{-1/3}x}
\left[\Phi(\beta_0,s)-\Phi(\beta_0,x_0) \right]\,\mathrm{d}s\\
&=
n^{2/3} \int_{x_0}^{x_0+n^{-1/3}x}
\left[ \Phi(\beta_0,s)-\Phi(\beta_0,x_0) \right]
\left[\lambda_0(s)-\lambda_0(x_0) \right]\,\mathrm{d}s\\
&=
\O_p(n^{-1/3}x),
\end{split}
\end{equation}
by conditions~\eqref{eq:cond Phi} and~\eqref{eq:cond lambda0}.
We conclude that
\begin{equation}
\label{eq:bound diff Zn LS ML}
\Phi(\beta_0,x_0) \left| \widetilde{\mathbbm{Z}}_n^{\lambda}(x)-\widehat{\mathbbm{Z}}_n^{\lambda}(x)  \right|=\O_p(n^{-1/2}x^{1/2})+\O_p(n^{-1/6}x)+\O_p(n^{-1/3}),
\end{equation}
and after taking the supremum over $[0,k]$, the lemma follows.
\end{proof}
To find the limit process of $\widehat{\mathbbm{Z}}_n^{\lambda}$,
we will apply results
from~\cite{kimpollard:1990}.
The limit distribution for~$\widetilde{\mathbbm{Z}}_n^{\lambda}$ will then follow directly from Lemma~\ref{lemma:relation_Z_n_tilde_hat}.
Let $\textbf{B}_{loc}(\mathbb{R})$ be the space of all locally bounded real functions on $\mathbb{R}$, equipped with the topology of uniform convergence on compact domains.
\begin{lemma}
\label{lemma:conv_distr_hat_Z_n}
Assume (A1) and (A2) and let $0<x_0<\tau_H$.
Suppose that~\eqref{eq:cond H^uc}, \eqref{eq:cond lambda0} and~\eqref{eq:cond Phi} hold.
Then the processes $\widehat{\mathbbm{Z}}_n^{\lambda}$ and $\widetilde{\mathbbm{Z}}_n^{\lambda}$
defined in~\eqref{def:Zn incr haz MLE} and~\eqref{def:Zn incr haz LS}
converge in distribution to the process
\begin{equation}
\label{def:Z}
\mathbbm{Z}(x)=\mathbbm{W} \left(\frac{\lambda_0(x_0)}{\Phi(\beta_0,x_0)}x\right)+\frac{1}{2}\lambda_0'(x_0)x^2,
\end{equation}
in $\textbf{B}_{loc}(\mathbb{R})$, where $\mathbbm{W}$ is standard two-sided Brownian motion originating from zero.
\end{lemma}
\begin{proof}
We will apply Theorem~4.7 in~\cite{kimpollard:1990}.
To this end, write the process $\widehat{\mathbbm{Z}}_n^{\lambda}$
in~\eqref{def:Zn incr haz MLE} as
\begin{equation}
\label{eq:Pn repr Znhat incr haz}
\widehat{\mathbbm{Z}}_n^{\lambda}(x)
=
-n^{2/3}\mathbbm{P}_ng(\cdot,n^{-1/3}x)+ n^{2/3}R_n(x),
\end{equation}
for $x\in[-n^{1/3}(x_0-T_{(1)}),n^{1/3}(T_{(n)}-x_0)]$,
where for $Y=(T,\Delta,Z)$ and $\theta\in[-x_0,\tau_H-x_0]$,
\begin{equation}
\label{def:g MLE}
\begin{split}
g(Y,\theta)&=-g_1(Y,\theta)+g_2(Y,\theta),\\
g_1(Y,\theta)&=\left( \{T< x_0+\theta\}-\{T< x_0\}\right) \frac{ \Delta}{\Phi(\beta_0,x_0)}\\
g_2(Y,\theta)&= \frac{\lambda_0(x_0)\text{e}^{\beta_0'Z}}{\Phi(\beta_0,x_0)} \int_{x_0}^{x_0+\theta} \{T\geq s\}\,\mathrm{d}s.
\end{split}
\end{equation}
Furthermore,
\[
R_n(x)
=
\frac{-\lambda_0(x_0)}{\Phi(\beta_0,x_0)}
\left[\left(\widehat{W}_n(x_0+n^{-1/3}x)- W_{n0}(x_0+n^{-1/3}x)\right)\\
-
\left(\widehat{W}_n(x_0)-W_{n0}(x_0)\right) \right],
\]
where $W_{n0}(x)=W_n(\beta_0,x)$, with $W_n$ defined in~\eqref{def:W_n}.
For all $k=1,2,\ldots,$ consider
\[
|R_n(x)|\leq \frac{\lambda_0(x_0)}{\Phi(\beta_0,x_0)} \int \left| \{ s\leq x_0+n^{-1/3}x \}-\{ s\leq x_0 \}\right| \left| \Phi_n(\hb,s)-\Phi_n(\beta_0,s) \right|\,\mathrm{d}s,
\]
which by similar reasoning as in~\eqref{eq:bound A_n} gives that
\begin{equation}
\label{eq:conv Rn MLE}
|R_n(x)|=\O_p(n^{-5/6}x),
\end{equation}
by Lemma~\ref{lemma:bounds Phin}.
Hence, the process $x\mapsto n^{2/3}R_n(x)$ tends to zero in $\textbf{B}_{loc}(\mathbb{R})$.
It is sufficient then to demonstrate that~$-n^{2/3}\mathbbm{P}_ng(\cdot,n^{-1/3}x)$ converges to
$\mathbbm{Z}(x)$ in $\textbf{B}_{loc}(\mathbb{R})$.
To this end, we will show that the conditions of Lemma~4.5 and~4.6 in~\cite{kimpollard:1990} hold.
Condition (i) of Lemma 4.5 is trivially fulfilled, since $\theta_0=0$ is an interior point of $[-x_0,\tau_H-x_0]$.
Moreover, observe that for all $\theta\in[-x_0,\tau_H-x_0]$,
from~\eqref{eq:relation base haz} and~\eqref{def:Phi}, we have

\begin{equation}
\label{eq:Pg}
Pg(\cdot,\theta)
=
\frac{-1}{\Phi(\beta_0,x_0)}\int_{x_0}^{x_0+\theta}
\left[\lambda_0(u)-\lambda(x_0)\right]\Phi(\beta_0,u)\, \mathrm{d}u.
\end{equation}
Thus, by~\eqref{eq:cond Phi} and~\eqref{eq:cond lambda0},
\[
\begin{split}
\frac{\partial Pg(\cdot,\theta)}{\partial \theta}
&=
-\frac{\Phi(\beta_0,x_0+\theta)}{\Phi(\beta_0,x_0)}\left\{ \lambda_0(x_0+\theta)-\lambda_0(x_0) \right\}\\
\frac{\partial^2 Pg(\cdot,\theta)}{\partial \theta^2}
&=
-\left(\frac{\partial \Phi(\beta_0,x_0+\theta)}{\partial\theta}\right)
\frac{\lambda_0(x_0+\theta)-\lambda_0(x_0)}{\Phi(\beta_0,x_0)}
-
\frac{\Phi(\beta_0,x_0+\theta)}{\Phi(\beta_0,x_0)}\lambda'_0(x_0+\theta).
\end{split}
\]%
It follows that $Pg(\cdot,\theta)$ is twice differentiable at $\theta_0=0$,
its unique maximizing value, with second derivative~$-\lambda_0'(x_0)<0$,
which establishes condition (iii) of Lemma~4.5 in~\cite{kimpollard:1990}.
Next, compute
\[
H(s,t)=\lim_{\alpha\to\infty} \alpha Pg(\cdot,s/\alpha)g(\cdot,t/\alpha),
\]
for finite $s$ and $t$.
Write
\[
\alpha Pg(\cdot,s/\alpha)g(\cdot,t/\alpha)
=
\alpha P
\Big(-g_1(\cdot,s/\alpha)+g_2(\cdot,s/\alpha)\Big)
\Big(-g_1(\cdot,t/\alpha)+g_2(\cdot,t/\alpha)\Big)
\]
and compute the four terms separately.
For all $s$ and $t$,
\begin{equation}
\label{eq:lim Pg1g2}
\alpha P\left|g_1(\cdot,s/\alpha)g_2(\cdot,t/\alpha)\right|\leq \frac{\lambda_0(x_0)t}{\Phi^2(\beta_0,x_0)}
\mathbb{E}
\left[
|\{T<x_0+s/\alpha\}-\{T<x_0\}|
\text{e}^{\beta_0'Z}
\right]
\to0,
\end{equation}
as $\alpha\to\infty$.
Completely analogous, it follows that
\begin{equation}
\label{eq:lim Pg2g2}
\lim_{\alpha\to\infty} \alpha
Pg_2(\cdot,s/\alpha)g_2(\cdot,t/\alpha)=0,
\end{equation}
for all $s$ and $t$.
Finally, consider the limit for $\alpha P g_1(\cdot,s/\alpha)g_1(\cdot,t/\alpha)$.
For $s,t\geq0$,
\[
\begin{split}
\alpha P g_1(\cdot,s/\alpha)g_1(\cdot,t/\alpha)
&=
\frac{\alpha}{\Phi^2(\beta_0,x_0)}\int \delta\{ x_0\leq u<x_0+(s\wedge t)/\alpha\}\, \mathrm{d} P(u,\delta,z)\\
&=
\frac{\alpha}{\Phi^2(\beta_0,x_0)} \int_{x_0}^{x_0+(s\wedge t)/\alpha} \lambda_0(u) \Phi(\beta_0,u)\, \mathrm{d}u\\
&=
\frac{1}{\Phi^2(\beta_0,x_0)} \int_0^{s\wedge t} \lambda_0(x_0+v/\alpha) \Phi(\beta_0,x_0+v/\alpha)\, \mathrm{d}v,
\end{split}
\]
by~\eqref{eq:relation base haz} and~\eqref{def:Phi}.
Therefore, by the continuity of $\lambda_0$ and $\Phi$,
\begin{equation}
\label{eq:lim Pg1g1}
\lim_{\alpha\to\infty} \alpha Pg_1(\cdot,s/\alpha)g_1(\cdot,t/\alpha)=\frac{\lambda_0(x_0)}{\Phi(\beta_0,x_0)} (s\wedge t).
\end{equation}
A similar reasoning applies for $s,t<0$
and $P g_1(\cdot,s/\alpha)g_1(\cdot,t/\alpha)=0$, when $s$ and $t$ have opposite signs.
Hence, condition (ii) of Lemma~4.5 in~\cite{kimpollard:1990} is verified, with
\[
H(s,t)=\frac{\lambda_0(x_0)}{\Phi(\beta_0,x_0)} (|s|\wedge |t|),
\]
for $st\geq0$ and $H(s,t)=0$, for $st<0$.
Note that $H(s,t)$ is the covariance kernel of the centered Gaussian process in~\eqref{def:Z}.
For condition (iv) of Lemma~4.5 in~\cite{kimpollard:1990}, it needs to be shown that for each $t$ and $\varepsilon>0$
\begin{equation}
\label{eq:L45iv_MLE}
\lim_{\alpha\to\infty}\alpha
Pg(\cdot,t/\alpha)^2\{|g(\cdot,t/\alpha)|>\alpha\varepsilon\}=0.
\end{equation}
In view of \eqref{eq:lim Pg1g2} and \eqref{eq:lim Pg2g2}, it suffices to show that
\[
\lim_{\alpha\to\infty}\alpha
Pg_1(\cdot,t/\alpha)^2\{|g(\cdot,t/\alpha)|>\alpha\varepsilon\}=0.
\]
Moreover, since $g_1$ is bounded uniformly for $\theta\in[-x_0,\tau_H-x_0]$,
by Lemma~\ref{lemma:bounds Phin},
\[
\{|g(\cdot,t/\alpha)|>\alpha \varepsilon \}\leq
\{|g_2(\cdot,t/\alpha)|>\alpha \varepsilon/2 \}\leq
\frac{2}{\alpha \varepsilon}|g_2(\cdot,t/\alpha)|,
\]
for $\alpha$ sufficiently large.
By \eqref{eq:lim Pg1g2}, it follows that
\[\begin{split}
\alpha Pg_1(\cdot,t/\alpha)^2 \{|g(\cdot,t/\alpha)|>\alpha \varepsilon \} &\leq  \frac{2}{\varepsilon}Pg_1(\cdot,t/\alpha)^2\left| g_2(\cdot,t/\alpha) \right|\\
&\leq \frac{2}{\varepsilon\Phi(\beta_0,M)}P\left| g_1(\cdot,t/\alpha)g_2(\cdot,t/\alpha) \right|\to 0.\\
\end{split}\]
Hence all conditions of Lemma~4.5 in~\cite{kimpollard:1990} are satisfied.

To continue with verifying the conditions of Lemma~4.6 in~\cite{kimpollard:1990},
consider the class of functions $\mathcal{G}=\{g(\cdot,\theta): \theta\in[-x_0,\tau_H-x_0]\}$ and the classes
\begin{equation}
\label{eq:class GR MLE}
\mathcal{G}_R=\left\{ g(\cdot,\theta)\in \mathcal{G}: |\theta|\leq R \right\},
\end{equation}
for any $R>0$, $R$ in a neighborhood of zero.
Since the functions in $\mathcal{G}_R$ are the
difference of~$g_1(\cdot,\theta)$, which is an the product of indicators, and $g_2(\cdot,\theta)$,
which is the product of a fixed function and a linear function,
it follows that $\mathcal{G}_R$ is a VC-subgraph class of functions, and hence it is uniformly manageable, which proves
condition (i) of Lemma~4.6 in~\cite{kimpollard:1990}.
Furthermore, choose as an envelope for $\mathcal{G}_R$,
\begin{equation}
\label{eq:envelope GR}
G_R=G_{R1}+G_{R2},
\end{equation}
where
\begin{equation}
\label{eq:def GR12_MLE}
\begin{split}
G_{R1}(T,\Delta,Z)
&=
\frac{\{x_0-R\leq T< x_0+R\}}{\Phi(\beta_0,x_0)},\\
G_{R2}(T,\Delta,Z)
&=
\frac{2R\lambda_0(x_0)}{\Phi(\beta_0,x_0)} \text{e}^{\beta_0'Z}.
\end{split}
\end{equation}%
Calculations completely analogous to \eqref{eq:lim Pg1g2}, \eqref{eq:lim Pg2g2} and \eqref{eq:lim Pg1g1}, with $1/R$ playing the role of $\alpha\to\infty$, yield that $PG^2_R=\O(R)$, as $R\to0$. This proves condition (ii) of Lemma~4.6 in~\cite{kimpollard:1990}. %which shows that $G_R$ in \eqref{eq:envelope GR} is square integrable, for $R$ in a neighborhood of $0$.
To show condition (iii) of Lemma~4.6 in~\cite{kimpollard:1990}, first note that
\[
P|g(\cdot,\theta_1)-g(\cdot,\theta_2)|\leq P|g_1(\cdot,\theta_1)-g_1(\cdot,\theta_2)|+P|g_2(\cdot,\theta_1)-g_2(\cdot,\theta_2)|.
\]
Now,
\[
P|g_1(\cdot,\theta_1)-g_1(\cdot,\theta_2)|
=
\frac{1}{\Phi(\beta_0,x_0)}\left| H^{uc}(x_0+\theta_1)-H^{uc}(x_0+\theta_2) \right|
=
\O(|\theta_1-\theta_2|),
\]
according to~\eqref{eq:cond H^uc}.
Analogously,
\[
P|g_2(\cdot,\theta_1)-g_2(\cdot,\theta_2)|
\leq
\frac{\lambda_0(x_0)}{\Phi(\beta_0,x_0)} |\theta_1-\theta_2|
\mathbb{E}\left[\text{e}^{\beta_0'Z}\right]
=\O(|\theta_1-\theta_2|),
\]
by~(A2), which proves condition (iii) of Lemma~4.6 in~\cite{kimpollard:1990}.
Finally, to establish condition (iv) of Lemma~4.6 in~\cite{kimpollard:1990}, we have to show that for each $\varepsilon>0$, there exists $K>0$ such that
\[
PG_R^2\{G_R>K\}<\varepsilon R,
\]
for $R$ near zero. The proof of this is completely analogous to proving \eqref{eq:L45iv_MLE}, with $1/R$ playing the role $\alpha\to\infty$. This shows that all conditions of Theorem~4.7 in \cite{kimpollard:1990} are fulfilled, from which we conclude that the process $-n^{2/3}\mathbbm{P}_ng(\cdot,n^{-1/3}x)$ converges in distribution to the process
\[
-\mathbbm{W}\left(\frac{\lambda_0(x_0)}{\Phi(\beta_0,x_0)}x\right)+\frac{1}{2}\lambda_0'(x_0)x^2
%\stackrel{d}{=}
\eqd
\mathbbm{W}\left(\frac{\lambda_0(x_0)}{\Phi(\beta_0,x_0)}x\right)+\frac{1}{2}\lambda_0'(x_0)x^2.
\]
Together with \eqref{eq:Pn repr Znhat incr haz} and~\eqref{eq:conv Rn MLE},
this proves the weak convergence of~$\widehat{\mathbbm{Z}}_n^{\lambda}$.
Weak convergence of~$\widetilde{\mathbbm{Z}}_n^{\lambda}$ is then immediate, by Lemma~\ref{lemma:relation_Z_n_tilde_hat}.
\end{proof}
As a consequence, we obtain the limiting distribution of the process in~\eqref{def:Un scaled haz ML}.
\begin{lemma}
\label{lemma:conv_distr_whole_proc}
Assume (A1) and (A2) and suppose that~\eqref{eq:cond H^uc}, \eqref{eq:cond lambda0} and~\eqref{eq:cond Phi} hold.
Let $0<x_0<\tau_H$ and $a>0$ fixed and let $\widehat{\mathbbm{Z}}_n^{\lambda}$ and $\widehat{W}_n$
be defined in~\eqref{def:Zn incr haz MLE} and~\eqref{def:Wn_hat}.
Then, the process
\[
\widehat{\mathbbm{Z}}_n^{\lambda}(x)-\frac{n^{1/3}a}{\Phi(\beta_0,x_0)}
\left[ \widehat{W}_n(x_0+n^{-1/3}x)-\widehat{W}_n(x_0) \right]
\]
converges weakly to
\[
\mathbbm{Z}(x)-ax=\mathbbm{W}\left(\frac{\lambda_0(x_0)}{\Phi(\beta_0,x_0)}x\right)+\frac{1}{2}\lambda_0'(x_0)x^2 -ax,
\]
in $\textbf{B}_{loc}(\mathbb{R})$, where $\mathbbm{W}$ is standard two-sided Brownian motion originating from zero.
\end{lemma}
\begin{proof}
In view of Lemma~\ref{lemma:conv_distr_hat_Z_n}, it suffices to show that for any $k=1,2,\ldots$,
\begin{equation}
\label{eq:drift Znhat}
\sup_{|x|\leq k}
\left|
n^{1/3}
\left[
\widehat{W}_n(x_0+n^{-1/3}x)-\widehat{W}_n(x_0)
\right] - \Phi(\beta_0,x_0) x \right|\to 0,
\end{equation}
almost surely.
This is immediate, since similar to~\eqref{eq:bound Wnhat-W0 main text},
together with the monotonicity of~$\Phi(\beta_0,u)$, one has, for $x\geq0$,
\begin{equation}
\label{eq:bound diff Wnhat x}
\begin{split}
\Big| n^{1/3}
&
\left[
\widehat{W}_n(x_0+n^{-1/3}x)-\widehat{W}_n(x_0)
\right] -\Phi(\beta_0,x_0) x \Big|\\
&\leq
n^{1/3}
\int_{x_0}^{x_0+n^{-1/3}x}
\left|
\Phi_n(\hb,u)-\Phi(\beta_0,x_0)
\right|
\,\mathrm{d}u\\
&\leq
|x| \sup_{u\in\mathbbm{R}}\left| \Phi_n(\hb,u)-\Phi(\beta_0,u) \right|
+
|\Phi(\beta_0,x_0+n^{-1/3}x)-\Phi(\beta_0,x_0)|\\
&=
o(x)+\O(n^{-1/3}x),
\end{split}
\end{equation}
almost surely, using~Lemma~\ref{lemma:conv Phin} and~\eqref{eq:cond Phi}.
The case $x<0$ can be treated likewise.
\end{proof}

Finally, the next lemma provides the limit process of $\widetilde{\mathbbm{Z}}_n^f$.
\begin{lemma}
Assume (A1) and (A2). Let $x_0\in(0,\tau_H)$ and suppose that~\eqref{eq:cond H^uc},
\eqref{eq:cond lambda0} and~\eqref{eq:cond Phi} hold.
Then the process $\widetilde{\mathbbm{Z}}_n^f$ defined in~\eqref{def:Zn dens} converges in distribution to the process
\begin{equation}
\label{def:Zf}
\mathbbm{Z}^f(x)= \mathbbm{W}\left( \frac{f_0(x_0)(1-F_0(x_0))}{\Phi(\beta_0,x_0)}x \right) +\frac{1}{2}f_0'(x_0)x^2.
\end{equation}
in $\textbf{B}_{loc}(\mathbb{R})$, where $\mathbbm{W}$ is standard two-sided Brownian motion originating from zero.
\end{lemma}
\begin{proof}
From~\eqref{def:Zn incr haz LS}, we have
$\BR(x_0+n^{-1/3}x)-\BR(x_0)=n^{-2/3}\widetilde{\mathbbm{Z}}_n^{\lambda}(x)+n^{-1/3}\lambda_0(x_0) x$,
so that by~\eqref{eq: relation Fn},
\begin{equation}
\label{eq:Znf in Znlambda}
\begin{split}
\mathbbm{\widetilde{Z}}^f_n(x)
&=
n^{2/3}
\left[-\text{e}^{-\BR(x_0+n^{-1/3}x)}+\text{e}^{-\BR(x_0)}
-n^{-1/3}f_0(x_0) x\right]\\
&=
n^{2/3}
\left[-\text{e}^{-\BR(x_0)}
\left(
\text{e}^{-n^{-2/3}\widetilde{\mathbbm{Z}}_n^{\lambda}(x)-n^{-1/3}\lambda_0(x_0) x}-1\right)
-n^{-1/3}f_0(x_0) x\right].
\end{split}
\end{equation}
Because $\text{e}^{-y}-1=-y+y^2/2+o(y^2)$, for $y\to0$ and $\sup_{x\in\mathbb{R}}|\widetilde{\mathbbm{Z}}_n^{\lambda}(x)|=\O_p(1)$,
according to~Lemma~\ref{lemma:conv_distr_hat_Z_n},
it follows that
\[\begin{split}
\text{e}^{-n^{-2/3}\widetilde{\mathbbm{Z}}_n^{\lambda}(x)-n^{-1/3}\lambda_0(x_0) x}-1
=
-n^{-2/3}\widetilde{\mathbbm{Z}}_n^{\lambda}(x)
&-
n^{-1/3}\lambda_0(x_0) x
+
\frac{1}{2}n^{-2/3}\lambda_0(x_0)^2 x^2\\
&+
\O_p(n^{-4/3})+\O_p(n^{-1}x)+o_p(n^{-2/3}x^2).
\end{split}
\]
Similarly, from Theorem~\ref{theorem:conv_breslow}, we have
that $\text{e}^{-\BR(x_0)}=\text{e}^{-\Lambda_0(x_0)}+\O_p(n^{-1/2})$.
Since
\[
\text{e}^{-\Lambda_0(x_0)}\lambda_0(x_0)=(1-F_0(x_0))\lambda_0(x_0)=f_0(x_0),
\]
from~\eqref{eq:Znf in Znlambda}, we find that
\begin{equation}
\label{eq:Znf in Znlambda expansion}
\begin{split}
\mathbbm{\widetilde{Z}}^f_n(x)
=
(1-F_0(x_0))\widetilde{\mathbbm{Z}}_n^{\lambda}(x)
&-
\frac{1}{2}(1-F_0(x_0))\lambda_0(x_0)^2 x^2\\
&+
\O_p(n^{-1/2})
+
\O_p(n^{-1/6}x)
+
o_p(x^2).
\end{split}
\end{equation}
According to Lemma~\ref{lemma:conv_distr_hat_Z_n}, the process
$(1-F_0(x_0))\widetilde{\mathbbm{Z}}_n^{\lambda}(x)
-
\frac{1}{2}(1-F_0(x_0))\lambda_0(x_0)^2 x^2$
converges weakly to
\[
(1-F_0(x_0))\mathbbm{W}\left( \frac{\lambda_0(x_0)}{\Phi(\beta_0,x_0)}x \right)
+
\frac{1}{2}(1-F_0(x_0))\lambda_0'(x_0)x^2-\frac{1}{2}(1-F_0(x_0))\lambda_0^2(x_0)x^2,
\]
which has the same distribution as the process in~\eqref{def:Zf},
by means of Brownian scaling and the fact that
\begin{equation}
\label{eq:relation lambda' and f'}
\lambda_0'=\left(\frac{f_0}{1-F_0}\right)'=\frac{(1-F_0)f_0'+f_0^2}{(1-F_0)^2}
=
\frac{f_0'}{1-F_0}+\lambda_0^2.
\end{equation}
Hence, for any $k=1,2,\ldots$, it follows from~\eqref{eq:Znf in Znlambda expansion} that
\[
\sup_{|x|\leq k}|\mathbbm{\widetilde{Z}}^f_n(x)-\mathbbm{Z}^f(x)|=o_p(1),
\]
which finishes the proof.
\end{proof}

\section{Limit distribution}
\label{sec:limit distribution}
\noindent
The last step in deriving the asymptotic distribution of the estimators is to find the limiting distribution of the
inverse processes $\widetilde{U}_n^{\lambda}$, $\widehat{U}_n^{\lambda}$ and $\widetilde{U}_n^f$
defined in~\eqref{def:inverse LS haz}, \eqref{def:inverse MLE} and~\eqref{def:inverse_process_f}
and of the versions of $\widetilde{U}_n^{\lambda}$ and $\widehat{U}_n^{\lambda}$ in the case of a nonincreasing hazard,
by applying Theorem~2.7 in~\cite{kimpollard:1990}.
This requires the inverse processes to be bounded in probability.
\begin{lemma}
\label{lemma:bound_local_proc}
Assume (A1) and (A2) and let $x_0\in(0,\tau_H)$.
Suppose that $\lambda_0$ is monotone and suppose that $f_0$ is nondecreasing.
Suppose that~\eqref{eq:cond lambda0} and~\eqref{eq:cond Phi} hold,
with $\lambda_0(x_0)\ne 0$.
Then, for each $\varepsilon>0$ and $M_1>0$, there exists $M_2>0$ such that, for $n$ large enough,
\begin{align}
\mathbbm{P}\left( \max_{|a|\leq M_1} n^{1/3}\left|\widehat{U}_n^{\lambda}(\lambda_0(x_0)+n^{-1/3}a)-x_0\right|>M_2 \right)
&<\varepsilon
\label{eq:bounded Uhat lambda}\\
\mathbbm{P}\left( \max_{|a|\leq M_1} n^{1/3}\left|\widetilde{U}_n^{\lambda}(\lambda_0(x_0)+n^{-1/3}a)-x_0\right|>M_2 \right)
&<\varepsilon
\label{eq:bounded Utilde lambda}\\
\mathbbm{P}\left( \max_{|a|\leq M_1} n^{1/3}\left|\widetilde{U}_n^f(f_0(x_0)+n^{-1/3}a)-x_0\right|>M_2 \right)
&<\varepsilon
\label{eq:bounded Utilde f},
\end{align}
for $n$ sufficiently large.
\end{lemma}
\noindent
The proof can be found in the appendix.
Hereafter, the continuous mapping theorem from~\cite{kimpollard:1990}
will be applied to the inverse processes in~~\eqref{def:inverse MLE}, \eqref{def:inverse LS haz} and~\eqref{def:inverse_process_f},
in order to derive the limiting distribution of the considered estimators.
Let $\mathbbm{C}_{max}(\mathbb{R})$ denote the subset of $\textbf{B}_{loc}(\mathbb{R})$
consisting of continuous functions $f$ for which $f(t)\to -\infty$, when $|t|\to\infty$ and $f$ has an unique maximum.

\begin{proof}[ Proof of Theorem~\ref{theorem:asymp_hazard_nondecr}]
The aim is to apply Theorem~2.7 in~\cite{kimpollard:1990} and Theorem~6.1 in~\cite{huangwellner:1995}.
Since Theorem~2.7 from~\cite{kimpollard:1990} applies to the argmax of processes on the whole real line,
we extend the process
\[
\widehat{Z}_n^\lambda(a,x)
=
 \widehat{\mathbbm{Z}}_n^{\lambda}(x)- \frac{n^{1/3}a}{\Phi(\beta_0,x_0)} \left[ \widehat{W}_n(x_0+n^{-1/3}x)-\widehat{W}_n(x_0) \right]
\]
from~\eqref{def:Un scaled haz ML} for $x\in [n^{1/3}(T_{(1)}-x_0),n^{1/3}(T_{(n)}-x_0)]$, to the whole real line.
Define $\widehat{Z}_n^{\lambda}(a,x)=\widehat{Z}_n^{\lambda}(a,n^{1/3}(T_{(1)}-x_0))$,
for $x<n^{1/3}(T_{(1)}-x_0)$ and
$\widehat{Z}_n^{\lambda}(a,x)=\widehat{Z}_n^{\lambda}(a,n^{1/3}(T_{(n)}-x_0))+1$,
for $x>n^{1/3}(T_{(n)}-x_0)$.
Then, $\widehat{Z}_n^{\lambda}(a,x)\in \textbf{B}_{loc}(\mathbbm{R})$ and according to~\eqref{def:Un scaled haz ML},
\[
n^{1/3}
\left[\widehat{U}_n^{\lambda}(\lambda_0(x_0)+n^{-1/3}a)-x_0\right]
=
\argmin_{x\in\mathbb{R}}
\left\{
\widehat{Z}_n^{\lambda}(a,x)
\right\}
=
\argmax_{x\in \mathbbm{R}}
\left\{
-\widehat{Z}_n^{\lambda}(a,x)
\right\}.
\]
By Lemma~\ref{lemma:conv_distr_hat_Z_n}, for any $a$ fixed,
the process $-\widehat{Z}_n^{\lambda}(a,x)$
converges weakly to the process $-\mathbbm{Z}(x)+ax\in\mathbbm{C}_{max}(\mathbbm{R})$
with probability one,
where $\mathbbm{Z}$ has been defined in~\eqref{def:Z}.
Lemma~\ref{lemma:bound_local_proc} ensures the boundedness in probability of
$n^{1/3}\{\widehat{U}_n^{\lambda}(\lambda_0(x_0)+n^{-1/3}a)-x_0\}$.
Consequently, by Theorem~2.7 in~\cite{kimpollard:1990} it follows that
\[
n^{1/3}
\left[\widehat{U}_n^{\lambda}(\lambda_0(x_0)+n^{-1/3}a)-x_0\right]
\law
\argmax_{x\in\mathbb{R}}
\left\{
-\mathbbm{Z}(x)+ax
\right\}
=
\argmin_{x\in\mathbbm{R}}
\left\{
\mathbbm{Z}(x)-ax
\right\}.
\]
The same argument applies to the process
$\widetilde{\mathbbm{Z}}_n^{\lambda}(x)-ax$
from~\eqref{def:Un scaled haz LS}, for $x\in [-n^{1/3}x_0,n^{1/3}(T_{(n)}-x_0)]$,
which we extend to the whole real line in a similar fashion.
Furthermore, if we fix $a,b\in\mathbbm{R}$, it will follow that
\[
\left(
\widehat{Z}_n^{\lambda}(a,x), \widetilde{\mathbbm{Z}}_n^{\lambda}(x)-bx
\right)
\law
\Big(
\mathbbm{Z}(x)-ax,\mathbbm{Z}(x)-bx
\Big),
\]
by Lemma~\ref{lemma:conv_distr_whole_proc} and Lemma~\ref{lemma:conv_distr_hat_Z_n}.
Hence, the first condition of Theorem~6.1 in~\cite{huangwellner:1995} is verified.
The second condition is provided by Lemma~\ref{lemma:bound_local_proc}, whereas the third condition is given
by~\eqref{def:Un scaled haz LS} and~\eqref{def:Un scaled haz ML}.
Therefore, by Theorem~6.1 in~\cite{huangwellner:1995},
\[
\left( n^{1/3}
\left[\widehat{U}_n^{\lambda}(\lambda_0(x_0)+n^{-1/3}a)-x_0\right],
n^{1/3}\left[\widetilde{U}_n^{\lambda}(\lambda_0(x_0)+n^{-1/3}b)-x_0\right] \right)
\law
\left(U^{\lambda}(a),U^{\lambda}(b) \right),
\]
where
\[
U^{\lambda}(a)=\sup \left\{ t: \mathbbm{W}\left(\frac{\lambda_0(x_0)}{\Phi(\beta_0,x_0)}\,t\right)+\frac{1}{2}\lambda_0'(x_0)t^2 -at \text{ is minimal} \right\}.
\]
Additional computations show that $U^{\lambda}(a)\eqd U^{\lambda}(0)+a/\lambda_0'(x_0)$ and therefore,
by the definition of the inverse processes in \eqref{def:inverse MLE} and~\eqref{def:inverse LS haz},
\[
\begin{split}
&
\mathbbm{P}\left(n^{1/3}\left[\hoMLE(x_0)-\lambda_0(x_0)\right]>a,
n^{1/3}\left[\ho(x_0)-\lambda_0(x_0)\right]>b\right)\\
&\quad\to
\mathbbm{P}(U^{\lambda}(a)<0,U^{\lambda}(b)<0)
=
\mathbbm{P}(-\lambda_0'(x_0)U^{\lambda}(0)>a,-\lambda_0'(x_0)U^{\lambda}(0)>b),
\end{split}
\]
as $n\to \infty$.
This implies that
\[
\left(n^{1/3}\left[\hoMLE(x_0)-\lambda_0(x_0)\right],
n^{1/3}\left[\ho(x_0)-\lambda_0(x_0)\right]  \right)\law \left(-\lambda_0'(x_0)U^{\lambda}(0),-\lambda_0'(x_0)U^{\lambda}(0) \right),
\]
which proves~\eqref{eq:asymp equiv MLE-LS}.
To establish the limiting distribution,
define
\[
A(x)=\left( \frac{\Phi(\beta_0,x)}{4\lambda_0(x)\lambda_0'(x)} \right)^{1/3},
\]
and note that
\[
n^{1/3}A(x_0)\left[\hoMLE(x_0)-\lambda_0(x_0)\right]
\law
A(x_0)\lambda'_0(x_0)U^{\lambda}(0)
\stackrel{d}{=}
\argmin_{t\in \mathbbm{R}}\left\{ \mathbbm{W}(t)+t^2 \right\},
\]
by Brownian scaling and the fact that the distribution of $U^{\lambda}(0)$ is symmetric around zero.
\end{proof}
\begin{proof}[Proof of Theorem~\ref{theorem:asymp_hazard_nonincr}]
The proof of Theorem~\ref{theorem:asymp_hazard_nonincr} is completely analogous
to that of Theorem~\ref{theorem:asymp_hazard_nondecr}.
The inverse processes to be considered in this case are
\[
\begin{split}
\widehat{U}_n^{\lambda}(a)
&=
\argmax_{x\in[0,T_{(n)}]}\left\{Y_n(x)-aW_n(\hb,x)\right\},\\
\widetilde{U}_n^{\lambda}(a)
&=
\argmax_{x\in[0,T_{(n)}]}
\left\{\BR(x)-ax\right\},
\end{split}
\]
for $a>0$, where $W_n$, $Y_n$ and $\Lambda_n$ have been defined in~\eqref{def:W_n},
\eqref{def:Y_n} and~\eqref{def:breslow} and $\hb$ is the maximum partial likelihood estimator.
By the same arguments as used in the proof of Theorem \ref{theorem:asymp_hazard_nondecr},
the limiting distribution is expressed in terms of
\[\begin{split}
\underset{t\in \mathbbm{R}}{\operatorname{argmax}}\left\{ \mathbbm{W}(t)-t^2 \right\}
\eqd
\underset{t\in \mathbbm{R}}{\operatorname{argmax}} \left\{ -\mathbbm{W}(t)-t^2 \right\}
=
\underset{t\in \mathbbm{R}}{\operatorname{argmin}}\{\mathbbm{W}(t)+t^2\},
\end{split}\]
by properties of Brownian motion.
\end{proof}
\begin{proof}[Proof of Theorem~\ref{theorem:asymp_density}]
Completely similar to the reasoning in the proof of Theorem~\ref{theorem:asymp_hazard_nondecr},
we obtain
\[
n^{1/3}
\left[\widetilde{U}_n^f(f_0(x_0)+n^{-1/3}a)-x_0\right]
\law
U^f(a),
\]
where
\[
U^f(a)
=
\sup\left\{ t: \mathbbm{W}\left( \frac{f_0(x_0)(1-F_0(x_0))}{\Phi(\beta_0,x_0)}t \right)+\frac{1}{2}f_0'(x_0)t^2-at \text{ is maximal}\right\}.
\]
As before, by Brownian scaling,
$U^f(a)\eqd U^f(0)+a/f_0'(x_0)$ and together with~\eqref{eq:switching decr} we obtain
\[
\mathbbm{P}\left(n^{1/3}
\left[\fo(x_0)-f_0(x_0)\right]<a\right)
\to
\mathbbm{P}\left(-f_0'(x_0)U^f(0)<a\right).
\]
Similar to the proof of Theorem~\ref{theorem:asymp_hazard_nondecr}, with
\[
A(x)=\left| \frac{\Phi(\beta_0,x)}{4f_0(x)f_0'(x)(1-F_0(x))} \right|^{1/3},
\]
we conclude that $n^{1/3}A(x_0)[\fo(x_0)-f_0(x_0)]$ converges in distribution to
\[
A(x_0)f_0'(x_0)U^f(0)
=
\argmax_{t\in\mathbbm{R}}\{ \mathbbm{W}(t)-t^2 \}\eqd \argmin_{t\in\mathbbm{R}}\{ \mathbbm{W}(t)+t^2 \},
\]
using Brownian scaling and the fact that the distribution of $U^f(0)$ is symmetric around zero.
\end{proof}

\paragraph{Acknowledgements.}
We would like to thank the associate editor and two anonymous referees for their valuable comments and suggestions.
We would also like to thank Jon Wellner for his help with the proof of Lemma~\ref{lemma:Jon}.

\section*{Supporting information}
\noindent
Additional information for this article is available online in Appendix S1:
Proofs of Lemma~\ref{lemma:bounds Phin}, Lemma~\ref{lemma:conv Wn Vn hat} and Lemma~\ref{lemma:bound_local_proc}.

\noindent
Hendrik P. Lopuha\"a, Delft Institute of Applied Mathematics, 
Delft University of Technology,
Mekelweg 4, 2628CD, Delft, The Netherlands.\\
E-mail: h.p.lopuhaa@tudelft.nl

%\bibliographystyle{acm}
%\bibliography{cox_RL_TN}

\pagebreak

%\end{document}

\pagebreak
%\end{comment}

\centerline{\LARGE Shape constrained nonparametric estimators of the baseline}
\smallskip
\centerline{\LARGE distribution in Cox proportional hazards model}
\bigskip
\centerline{\LARGE Supplementary Material}
\bigskip
\centerline{Hendrik P.~Lopuha\"{a} and Gabriela F.~Nane}
\centerline{Delft University of Technology}
\bigskip
\centerline{\date{\today}}

\setcounter{equation}{71}

\section*{Appendix S1}
\begin{proof}[Proof of Lemma~\ref{lemma:bounds Phin}]
First, for every $x\leq M$ and $\beta\in\mathbb{R}^p$,
\begin{equation}
\label{eq:Phi inequality(1)}
0<\Phi(\beta,M)\leq \Phi(\beta,x)
\end{equation}
and for every $x\in \mathbb{R}$ and $|\beta-\beta_0|\leq \varepsilon$,
\begin{equation}
\label{eq:Phi inequality(2)}
\Phi(\beta,x)\leq \Phi(\beta,0)
\leq
\sup_{|\beta-\beta_0 |\leq \varepsilon} \mathbb{E}\left[\text{e}^{\beta'Z}\right]<\infty.
\end{equation}
Hence, by dominated convergence, for every $x\leq M$,
the function $\beta\mapsto\Phi(\beta,x)$ is continuous and therefore attains a minimum on the set $|\beta-\beta_0|\leq \varepsilon$.
Together with~\eqref{eq:Phi inequality(1)} and~\eqref{eq:Phi inequality(2)}, this proves~(i).

To show (ii), note that similar to~\eqref{eq:Phi inequality(1)} and~\eqref{eq:Phi inequality(2)},
for every $x\in [0,M]$ and $\beta\in \mathbb{R}^p$,
\begin{equation}
\label{eq:Phin inequality(1)}
\Phi_n(\beta,M)\leq \Phi_n(\beta,x)
\end{equation}
and for every $x\in\mathbb{R}$ and $\beta\in \mathbb{R}^p$,
\begin{equation}
\label{eq:Phin inequality(2)}
\Phi_n(\beta,x)\leq \Phi_n(\beta,0).
\end{equation}
Choose $\varepsilon>0$ from (A2) and let $\delta=\varepsilon/2\sqrt{p}$.
Strong consistency of $\beta^*_n$ yields that, for $n$ sufficiently large,
\[
\beta_{0j}-\delta\leq \beta^*_{nj}\leq \beta_{0j}+\delta,
\quad
\text{for all }j=1,2,\ldots,p,
\]
with probability one.
Next, consider all subsets $I_k=\{i_1,i_2,\ldots,i_k\}\subseteq\{1,2,\ldots,p\}=I$.
Then, for each $I_k$ fixed, on each event
\[
\bigcap_{j\in I_k}\{Z_{ij}\geq 0\}\bigcap_{l\in I\setminus I_k}\{Z_{il}<0\},
\quad
\text{where }
Z_i=(Z_{i1},\ldots,Z_{ip})'\in \mathbb{R}^p,
\]
we have
\[
\sum_{j\in I_k}(\beta_{0j}-\delta)Z_{ij}+\sum_{l\in I\setminus I_k}(\beta_{0j}+\delta)Z_{il}
\leq
\beta^{*'}_nZ
\leq
\sum_{j\in I_k}(\beta_{0j}+\delta)Z_{ij}+\sum_{l\in I\setminus I_k}(\beta_{0j}-\delta)Z_{il}.
\]
Define $\alpha_k,\gamma_k\in \mathbb{R}^p$ with coordinates
\[
\alpha_{kj}
=
\begin{cases}
\beta_{0j}-\delta, & j\in I_k,\\
\beta_{0j}+\delta, & j\in I\setminus I_k,
\end{cases}
\quad\text{and}\quad
\gamma_{kj}
=
\begin{cases}
\beta_{0j}+\delta, & j\in I_k,\\
\beta_{0j}-\delta, & j\in I\setminus I_k.
\end{cases}
\]
Then $|\beta_0-\alpha_k|\leq\varepsilon$ and $|\beta_0-\gamma_k|\leq\varepsilon$
and together with~\eqref{eq:Phin inequality(1)} and ~\eqref{eq:Phin inequality(2)}, we find that for every $x\leq M$,
\begin{equation}
\label{eq:Phin bound(1)}
\min_{I_k\subseteq I}
\left\{
\frac{1}{n}\sum_{i=1}^n
\{T_i\geq M\}
\text{e}^{\alpha_k'Z_i}
\right\}
\leq
\Phi_n(\beta^*_n,x)
\end{equation}
and for every $x\in\mathbb{R}$,
\begin{equation}
\label{eq:Phin bound(2)}
\Phi_n(\beta^*_n,x)
\leq
\max_{I_k\subseteq I}
\left\{
\frac{1}{n}\sum_{i=1}^n
\text{e}^{\gamma_k'Z_i}
\right\}.
\end{equation}
By (A2) and the law of large numbers,
\[
\min_{I_k\subseteq I}
\left\{
\frac{1}{n}\sum_{i=1}^n
\{T_i\geq M\}
\text{e}^{\alpha_k'Z_i}
\right\}
\to
\min_{I_k\subseteq I}
\mathbb{E}\left[\{T\geq M\}\text{e}^{\alpha_k'Z}\right]>0,
\]
with probability one and similarly,
\begin{equation}
\label{eq:Phin bound LLN}
\max_{I_k\subseteq I}
\left\{
\frac{1}{n}\sum_{i=1}^n
\text{e}^{\gamma_k'Z_i}
\right\}
\to
\max_{I_k\subseteq I}
\mathbb{E}\left[\text{e}^{\gamma_k'Z}\right]
\leq
\sup_{|\beta-\beta_0 |\leq \varepsilon} \mathbb{E}\left[\text{e}^{\beta'Z}\right]<\infty,
\end{equation}
with probability one.
This proves (ii).

To prove (iii),
it suffices to show that the inequalities hold componentwise.
For this, notice that for the $j$th element of the vector $D^{(1)}$,
\[
\sup_{x\in\mathbb{R}}
\sup_{|\beta-\beta_0| \leq  \varepsilon}
\left|\mathbbm{E} \left[ \{T\geq x\} Z_j\, \text{e}^{\beta'Z} \right]\right|
\leq
\sup_{|\beta-\beta_0|\leq\varepsilon} \mathbbm{E}\left[ |Z_j|\text{e}^{\beta'Z} \right] < \infty,
\]
by (A2).
Completely analogous, a similar inequality can be shown for each element of $D^{(2)}$.

Finally, to prove $(iv)$, note that similar to~\eqref{eq:Phin bound(2)} and~\eqref{eq:Phin bound LLN},
for the $j$th component of $D_n^{(1)}$, we can write
\[
\sup_{x\in\mathbb{R}}
\left|
D_{nj}^{(1)}(\beta^*_n,x)
\right|
\leq
\sum_{I_k\subseteq I}
\left[
\frac{1}{n}\sum_{i=1}^n
|Z_i|
\text{e}^{\gamma_k'Z_i}
\right]
\to
\mathbb{E}
\left[
|Z|
\text{e}^{\gamma_k'Z}
\right]
<\infty,
\]
with probability one, as $n$ tends to infinity.
Likewise, a similar result can be obtained for each element of $D_n^{(2)}$.
\end{proof}

\begin{proof}[Proof of Lemma~\ref{lemma:conv Wn Vn hat}]
By Glivenko-Cantelli,
\begin{equation}
\label{eq:conv V_n}
\sup_{x\in [T_{(1)},T_{(n)}]} \left| V_n(x)-V(x) \right|\to 0,
\end{equation}
almost surely, because of the continuity of $V$.
Furthermore,
\begin{equation}
\label{eq:Wnhat as int}
W_n(\hb,T_{(1)})
=
\int_0^{T_{(1)}}\Phi_n(\hb,s)\,\mathrm{d}s
=
T_{(1)}\Phi_n(\hb,T_{(1)})\to 0,
\end{equation}
almost surely,
since $\Phi_n(\hb,s)$ is bounded uniformly according to
Lemma~\ref{lemma:bounds Phin} and
$T_{(1)}\to 0$ with probability one, by the Borel-Cantelli lemma.
Moreover, if we define
\begin{equation}
\label{def:W_0}
W(\beta,x)
=
\int \left( \text{e}^{\beta' z} \int_0^x \{ u\geq s \}\,\mathrm{d}s \right)\, \mathrm{d}P(u,\delta,z),
\end{equation}
then we can write
\begin{equation}
\label{eq:relation W_0 and Phi}
W_0(x)
=
W(\beta_0,x)
=
\int_0^x \Phi(\beta_0,s)\,\text{d}s,
\end{equation}
where $\Phi$ is defined in~\eqref{def:Phi}.
It follows that
\begin{equation}
\label{eq:bound Wnhat-W0}
\begin{split}
\sup_{x\in[T_{(1)},T_{(n)}]} \left| \widehat{W}_n(x)-W_0(x) \right| &\leq \sup_{x \in[T_{(1)},T_{(n)}]} \left| \int_0^x \left( \Phi_n(\hb,s)-\Phi(\beta_0,s) \right)\,\mathrm{d}s \right|,\\
&\leq \tau_H \sup_{x \in\mathbbm{R}} \left| \Phi_n(\hb,x)-\Phi(\beta_0,x) \right| \to 0,
\end{split}
\end{equation}
with probability one, by Lemma~\ref{lemma:conv Phin}.

Take $\widehat{W}_n^{-1}$ to be the inverse of $\widehat{W}_n$,
which is well defined on $[0,\widehat{W}_n(T_{(n)})]$, since $\widehat{W}_n$ is strictly monotone on $[T_{(1)},T_{(n)}]$.
We first extend $\widehat{W}_n$ to $[T_{(1)},\infty)$ and $\widehat{W}_n^{-1}$ to $[0,\infty)$.
Define $\widehat{W}_n(t)=\widehat{W}_n(T_{(n)})+(t-T_{(n)})$, for all $t\geq T_{(n)}$, so that $\widehat{W}_n^{-1}(y)=T_{(n)}+(y-\widehat{W}_n(T_{(n)}))$, for $y\geq \widehat{W}_n(T_{(n)})$.
Similarly, take $W_0^{-1}$ to be the inverse of $W_0$, which is well-defined since~$W_0$ is strictly monotone
on $[0,\tau_H]$ and extend~$W_0$ and~$W_0^{-1}$ to~$[0,\infty)$, by defining
$W_0(t)=W_0(\tau_H)+(t-\tau_H)$, for all~$t\geq \tau_H$,
so that~$W_0^{-1}(y)=\tau_H+(y-W_0(\tau_H))$, for $y\geq W_0(\tau_H)$.
It follows that the extension $W_0^{-1}(y)$ is uniformly continuous on~$[0,\infty)$.
Immediate derivations give that
\begin{equation}
\label{eq:bound Wnhat-W0 inverse}
\sup_{0\leq y\leq \widehat{W}_n(T_{(n)})}\left| \widehat{W}_n^{-1}(y)-W_0^{-1}(y) \right|\to 0,
\end{equation}
with probability one. Furthermore, it can be inferred that
\[
\begin{split}
\delta_n
&=
\sup_{y\in [0,\widehat{W}_n(T_{(n)})]} \left| V_n\circ \widehat{W}_n^{-1}(y)-V\circ W_0^{-1}(y) \right|\\
&\leq
\sup_{y\in [0,\widehat{W}_n(T_{(n)})]} \left| (V_n-V)\circ \widehat{W}_n^{-1}(y) \right|+\sup_{y\in [0,\widehat{W}_n(T_{(n)})]}
\left| V\circ \widehat{W}_n^{-1}(y)-V\circ W_0^{-1}(y) \right|\\
&\leq
\sup_{t\in [T_{(1)},T_{(n)}]} \left| V_n(t)-V(t) \right|+\sup_{y\in [0,\widehat{W}_n(T_{(n)})]} \left| V\circ \left(\widehat{W}_n^{-1}(y)-W_0^{-1}(y)\right) \right|\\
&\to 0,
\end{split}
\]
almost surely, by~\eqref{eq:conv V_n}, \eqref{eq:bound Wnhat-W0 inverse} and the continuity of $V$.
According to~\eqref{eq:relation base haz} and~\eqref{eq:relation W_0 and Phi},
$\lambda_0$ can also be represented as
\begin{equation}
\label{eq:repr lambda0}
\lambda_0(x)=\frac{\mathrm{d} V(x)/\mathrm{d}x}{\mathrm{d}W_0(x)/\mathrm{d}x},
\end{equation}
which is well-defined for $x\in[0,\tau_H)$, since $\Phi$ is bounded away from zero, by Lemma~\ref{lemma:bounds Phin}.
Taking $x=W_0^{-1}(y)$, gives that
\[
\frac{\mathrm{d}V\left(W_0^{-1}(y)\right)}{\mathrm{d} y}=\lambda_0\left(W_0^{-1}(y)\right), \quad y\in[0,W_0(\tau_H)).
\]
Therefore, convexity of $\Lambda_0$ implies convexity of
$V\circ W_0^{-1}$ and subsequently of~$V\circ W_0^{-1}-\delta_n$.
Moreover, from the definition of $\delta_n$, it follows that for every $y\in[0,\widehat{W}_n(T_{(n)})]$,
\[
V\circ W_0^{-1}(y)-\delta_n \leq V_n\circ \widehat{W}_n^{-1}(y).
\]
As $\widehat{V}_n\circ \widehat{W}_n^{-1}(y)$ is the greatest convex function below $V_n\circ \widehat{W}_n^{-1}(y)$, we must have
\[
V\circ W_0^{-1}(y)-\delta_n\leq \widehat{V}_n\circ \widehat{W}_n^{-1}(y)\leq V_n\circ \widehat{W}_n^{-1}(y),
\]
for each $y\in[0,\widehat{W}_n(T_{(n)})]$. Re-writing the above inequalities leads to
\[
-\delta_n\leq \widehat{V}_n\circ \widehat{W}_n^{-1}(y)-V\circ W_0^{-1}(y)\leq V_n\circ \widehat{W}_n^{-1}(y)-V\circ W_0^{-1}(y)\leq \delta_n.
\]
Taking the supremum over $[0,\widehat{W}_n(T_{(n)})]$ then yields
\begin{equation}
\label{eq:bound Wnhat-W0 composition}
\sup_{y\in[0,\widehat{W}_n(T_{(n)})]} \left| \widehat{V}_n\circ \widehat{W}_n^{-1}(y)-V\circ W_0^{-1}(y) \right|\to 0,
\end{equation}
with probability one.
From~\eqref{eq:bound Wnhat-W0 inverse},~\eqref{eq:bound Wnhat-W0 composition} and the continuity of $V$,
we conclude that
\[\begin{split}
\sup_{t\in [T_{(1)},T_{(n)}]} \left| \widehat{V}_n(t)-V(t) \right|
&=
\sup_{y\in[0,\widehat{W}_n(T_{(n)})]} \left| \left(\widehat{V}_n-V\right)\circ \widehat{W}_n^{-1}(y) \right|\\
&\leq
\sup_{y\in[0,\widehat{W}_n(T_{(n)})]} \left| \widehat{V}_n\circ \widehat{W}_n^{-1}(y)-V\circ W_0^{-1}(y) \right|\\
&\qquad+\sup_{y\in[0,\widehat{W}_n(T_{(n)})]} \left| V\circ W_0^{-1}(y)-V\circ \widehat{W}_n^{-1}(y) \right|
\to 0,
\end{split}
\]
with probability one.
\end{proof}

\begin{proof}[Proof of Lemma~\ref{lemma:bound_local_proc}]
The proof of the lemma follows closely the lines of proof of Lemma~5.3 in~\cite{APPENDIXgroeneboomwellner}
(see also Lemma~7.1 in~\cite{APPENDIXhuangwellner:1995}).
First consider~\eqref{eq:bounded Uhat lambda} in case $\lambda_0$ is nondecreasing.
It will be shown that
\begin{equation}
\label{eq:prob1}
\mathbbm{P}\left( \max_{|a|\leq M_1} n^{1/3}
\left[\widehat{U}_n^{\lambda}(\lambda_0(x_0)+n^{-1/3}a)-x_0\right]>M_2 \right) <\varepsilon,
\end{equation}
as the other part can be proved similarly.
Because $\widehat{U}_n^{\lambda}(a)$ is nondecreasing,
the probability in~\eqref{eq:prob1} is equal to
\[
\mathbbm{P}\left( n^{1/3}
\left[\widehat{U}_n(\lambda_0(x_0)+n^{-1/3}M_1)-x_0\right]>M_2 \right).
\]
The relationship between the inverse process $\widehat{U}_n^{\lambda}$ and the process $\widehat{\mathbbm{Z}}_n^{\lambda}$ defined
in~\eqref{def:Zn incr haz MLE}, together with the fact that $\widehat{\mathbbm{Z}}_n^{\lambda}(0)=0$, implies that
\begin{equation}
\label{eq:Un to Zn MLE}
\begin{split}
&
\mathbbm{P}\left( n^{1/3}
\left[ \widehat{U}_n^{\lambda}( \lambda_0(x_0)+n^{-1/3}M_1)-x_0\right]>M_2 \right) \\
&\leq
\mathbbm{P}\left( \widehat{\mathbbm{Z}}_n^{\lambda}(x)-\frac{n^{1/3}M_1}{\Phi(\beta_0,x_0)}
\left[\widehat{W}_n(x_0+n^{-1/3}x)-\widehat{W}_n(x_0)\right] \leq 0, \text{ for some } x\geq M_2 \right).
\end{split}
\end{equation}
By condition \eqref{eq:cond lambda0}, there exists $M_0>0$ such that,
for any $x\in[T_{(1)},T_{(n)}]$ with $|x-x_0|\leq M_0$, $\lambda_0'(x)>0$ and $\lambda_0'(x)$ is close to $\lambda_0'(x_0)$. Take $n^{-1/3}x\leq M_0$.
From~\eqref{eq:Pn repr Znhat incr haz} and~\eqref{eq:drift Znhat},
\begin{equation}
\label{eq:Png MLE}
\begin{split}
&
\widehat{\mathbbm{Z}}_n^{\lambda}(x)
-
\frac{n^{1/3}M_1}{\Phi(\beta_0,x_0)}
\left[ \widehat{W}_n(x_0+n^{-1/3}x)-\widehat{W}_n(x_0) \right]\\
&=
-n^{2/3}\mathbbm{P}_ng(\cdot,n^{-1/3}x)-M_1x+ \widehat{R}_n(x),
\end{split}
\end{equation}
where
$\widehat{R}_n(x)=\O_p(n^{-1/6}x)+o(x)+\O(n^{-1/3}x)$,
by \eqref{eq:conv Rn MLE} and \eqref{eq:bound diff Wnhat x}.
Furthermore, for $0<R\leq M_0$, consider the class of functions $\mathcal{G}_R$ defined in \eqref{eq:class GR MLE}
along with its envelope defined in~\eqref{eq:envelope GR}.
It has been determined in the proof of Lemma~\ref{lemma:conv_distr_hat_Z_n}
that $\mathcal{G}_R$ is uniformly manageable for its envelope $G_R$ and that $PG_R^2=\O(R)$, for $0<R\leq M_0$.
Thus, Lemma~4.1 in~\cite{APPENDIXkimpollard:1990} states that for each $\delta>0$, there exist random variables $S_n=\O_p(1)$ such that
\begin{equation}
\label{eq:Lemma41 Kim Pollard}
|\mathbbm{P}_ng(\cdot,n^{-1/3}x)-Pg(\cdot,n^{-1/3}x)|\leq \delta n^{-2/3}x^2+n^{-2/3}S_n^2,
\end{equation}
for $n^{-1/3}x\leq M_0$.
Choose $\delta=\lambda_0'(x_0)/8$ in the above inequality.
It will result that
\[
-n^{2/3}(\mathbbm{P}_n-P)g(\cdot,n^{-1/3}x)
\geq
-\frac{1}{8}\lambda_0'(x_0)x^2-S_n^2.
\]
Furthermore, by~\eqref{eq:cond lambda0}, \eqref{eq:cond Phi} and~\eqref{eq:Pg},
\begin{equation}
\label{eq:Pg expansion}
\begin{split}
-n^{2/3}Pg(\cdot,n^{-1/3}x)
=
\frac{x^2}{2\Phi(\beta_0,x_0)}
\Bigg(
\lambda_0'(x_0&+\theta_n)\Phi(\beta_0,x_0+\theta_n)\\
&+
\left[
\lambda_0(x_0+\theta_n)-\lambda_0(x_0)
\right]
\Phi'(\beta_0,x_0+\theta_n)
\Bigg)
\end{split}
\end{equation}
for $|\theta_n|\leq n^{-1/3}x\leq M_0$, where $\Phi'(\beta_0,x)=\partial\Phi(\beta_0,x)/\partial x$.
From the choice of $M_0$ and since~$\lambda_0'(x_0)>0$,
we can find $K>0$ such that for any $x>K$,
\[
-n^{2/3}Pg(\cdot,n^{-1/3}x)
-
M_1x\geq\frac{1}{4}\lambda_0'(x_0)x^2,
\]
for $n$ sufficiently large.
We conclude that
\[\begin{split}
\widehat{\mathbbm{Z}}_n^{\lambda}(x)&-\frac{n^{1/3}M_1}{\Phi(\beta_0,x_0)}
\left[ \widehat{W}_n(x_0+n^{-1/3}x)-\widehat{W}_n(x_0) \right]\\
&=
-n^{2/3}\mathbbm{P}_ng(\cdot,n^{-1/3}x)-M_1x+\widehat{R}_n(t)\\
&=
-n^{2/3}(\mathbbm{P}_n-P)g(\cdot,n^{-1/3}x)-n^{2/3}Pg(\cdot,n^{-1/3}x)-M_1x+\widehat{R}_n(x)\\
&\geq
\frac{1}{8}\lambda_0'(x_0)x^2-S_n^2+\widehat{R}_n(x),
\end{split}\]
where $\widehat{R}_n(x)=\O_p(n^{-1/6}x)+o(x)+\O(n^{-1/3}x)$ and the
$\O_p$, $\O$ and $o$ terms do not depend on $x$.
It follows that for $x\geq M_2> K$,
\begin{equation}
\label{eq:lwbound}
\widehat{\mathbbm{Z}}_n^{\lambda}(x)-\frac{n^{1/3}M_1}{\Phi(\beta_0,x_0)}
\left[ \widehat{W}_n(x_0+n^{-1/3}x)-\widehat{W}_n(x_0) \right]
\geq
\frac{1}{8}\lambda_0'(x_0)x^2-S_n^2+o_P(1),
\end{equation}
where the $o_P$ term does not depend on $x$.
Then, $M_2$ can be chosen such that
\[
\mathbbm{P}\left( S_n^2\geq \frac{1}{8}\lambda_0'(x_0)M_2^2+o_P(1) \right) <\varepsilon,
\]
for $n$ sufficiently large.
We find that
\[\begin{split}
\mathbbm{P}&\left( \widehat{\mathbbm{Z}}_n^{\lambda}(x)
-
\frac{n^{1/3}M_1}{\Phi(\beta_0,x_0)}
\left[ \widehat{W}_n(x_0+n^{-1/3}x)-\widehat{W}_n(x_0) \right] \leq 0,
\text{ for some } M_2\leq x\leq n^{1/3}M_0\right) \\
&\leq
\mathbbm{P}\left( \frac{1}{8}\lambda_0'(x_0)x^2-S_n^2+o_P(1)\leq 0, \text{ for some } M_2\leq x\leq n^{1/3}M_0 \right)\\
&\leq
\mathbbm{P}\left( S_n^2\geq \frac{1}{8}\lambda_0'(x_0)x^2+o_P(1), \text{ for some } M_2\leq x\leq n^{1/3}M_0 \right)\leq\varepsilon,
\end{split}\]
for $n$ sufficiently large.

For $n^{-1/3}x>M_0$, we first show that
\begin{equation}
\begin{split}
\label{eq:inequality_MLE}
\widehat{\mathbbm{Z}}_n^{\lambda}(x)
&-
\frac{n^{1/3}M_1}{\Phi(\beta_0,x_0)}
\left[ \widehat{W}_n(x_0+n^{-1/3}x)-\widehat{W}_n(x_0) \right] \\
&\geq
\widehat{\mathbbm{Z}}_n(n^{1/3}M_0/2)
-
\frac{n^{1/3}M_1}{\Phi(\beta_0,x_0)}
\left[ \widehat{W}_n(x_0+M_0/2)-\widehat{W}_n(x_0) \right],
\end{split}
\end{equation}
with large probability, for $n$ sufficiently large.
Then,
\[
\mathbbm{P} \left( \widehat{\mathbbm{Z}}_n(n^{1/3}M_0/2)-
\frac{n^{1/3}M_1}{\Phi(\beta_0,x_0)}
\left[ \widehat{W}_n(x_0+M_0/2)-\widehat{W}_n(x_0) \right]
\leq 0
\right)
\]
can be bounded with the argument above.
Lemma~\ref{lemma:conv Wn Vn hat} and~\eqref{eq:conv V_n} yield that
$\widehat{V}_n(x_0+M_0/2)=V_n(x_0+M_0/2)+o(1)$, with probability one
and by definition
$V_n(x_0+n^{-1/3}x)\geq \widehat{V}_n(x_0+n^{-1/3}x)$, for all
$x_0+n^{-1/3}x\in[T_{(1)},T_{(n)}]$.
This implies that
\begin{equation}
\label{eq:convexity MLE}
\begin{split}
&
V_n(x_0+n^{-1/3}x)-V_n(x_0+M_0/2)\\
&\qquad\geq
\widehat{V}_n(x_0+n^{-1/3}x)-\widehat{V}_n(x_0+M_0/2)+o(1),\\
&\qquad\geq
\hoMLE(x_0+M_0/2)\left(\widehat{W}_n(x_0+n^{-1/3}x)-\widehat{W}_n(x_0+M_0/2)\right)+o(1),
\end{split}
\end{equation}
using the convexity of $\widehat{V}_n$.
To show~\eqref{eq:inequality_MLE},
note that by definition~\eqref{def:Zn incr haz MLE},
\[\begin{split}
&
\widehat{\mathbbm{Z}}_n^{\lambda}(x)
-
\frac{n^{1/3}M_1}{\Phi(\beta_0,x_0)}
\left[ \widehat{W}_n(x_0+n^{-1/3}x)-\widehat{W}_n(\hb,x_0) \right]\\
&\qquad-
\left[\widehat{\mathbbm{Z}}_n^{\lambda}(n^{1/3}M_0/2)-
\frac{n^{1/3}M_1}{\Phi(\beta_0,x_0)}
\left( \widehat{W}_n(x_0+M_0/2)-\widehat{W}_n(x_0) \right)\right]\\
&=
\frac{n^{2/3}}{\Phi(\beta_0,x_0)}
\Bigg[
V_n(x_0+n^{-1/3}x)-V_n(x_0+M_0/2)\\
&\qquad\qquad\qquad\qquad-
\left(\lambda_0(x_0)+n^{-1/3}M_1\right)
\left(\widehat{W}_n(x_0+n^{-1/3}x)-\widehat{W}_n(x_0+M_0/2)\right) \Bigg]\\
&\geq
\frac{n^{2/3}}{\Phi(\beta_0,x_0)}
\Bigg[
\left( \hoMLE(x_0+M_0/2)-\lambda_0(x_0)-n^{-1/3}M_1  \right)\\
&\qquad\qquad\qquad\qquad\times
\left(
\widehat{W}_n(x_0+n^{-1/3}x)-\widehat{W}_n(x_0+M_0/2)\right)
+o(1)
\Bigg]\\
&=
\frac{n^{2/3}}{\Phi(\beta_0,x_0)}
\Bigg[
\left( \lambda_0(x_0+M_0/2)-\lambda_0(x_0)-n^{-1/3}M_1 +o(1) \right)\\
&\qquad\qquad\qquad\qquad\times
\left(
W_0(x_0+n^{-1/3}x)-W_0(x_0+M_0/2)+o(1)\right)
+o(1)
\Bigg]>0,
\end{split}\]
for $n$ sufficiently large,
using~\eqref{eq:bound Wnhat-W0} and the fact that
$\lambda_0$ and $W_0$ are strictly increasing and $n^{-1/3}x>M_0$.
It follows that
\[\begin{split}
\mathbbm{P}&\left(
\widehat{\mathbbm{Z}}_n^{\lambda}(x)-\frac{n^{1/3}M_1}{\Phi(\beta_0,x_0)}
\left[
\widehat{W}_n(x_0+n^{-1/3}x)-\widehat{W}_n(x_0)
\right] \leq 0, \text{ for some }   x>n^{1/3}M_0 \right)\\
&\leq
\mathbbm{P}\left( \widehat{\mathbbm{Z}}_n^{\lambda}(n^{1/3}M_0/2)-
\frac{n^{1/3}M_1}{\Phi(\beta_0,x_0)}
\left[
\widehat{W}_n(x_0+M_0/2)-\widehat{W}_n(x_0)
\right] \leq 0 \right)
\leq \varepsilon.
\end{split}\]
This completes the proof of~\eqref{eq:prob1}.
The other part of~\eqref{eq:bounded Uhat lambda} for a nondecreasing $\lambda_0$ is proven similarly.

For~\eqref{eq:bounded Utilde lambda}, in case of a nondecreasing $\lambda_0$,
by the same reasoning that leads to~\eqref{eq:Un to Zn MLE} we first have
\[
\mathbbm{P}\left( n^{1/3}
\left[ \widetilde{U}_n^{\lambda}( \lambda_0(x_0)+n^{-1/3}M_1)-x_0\right]>M_2 \right)
\leq
\mathbbm{P}\left( \widetilde{\mathbbm{Z}}_n^{\lambda}(x)-M_1x \leq 0, \text{ for some } x\geq M_2 \right).
\]
Moreover, by~\eqref{eq:bound diff Zn LS ML},
\[
\widetilde{\mathbbm{Z}}_n^{\lambda}(x)=\widehat{\mathbbm{Z}}_n^{\lambda}(x)+\O_p(n^{-1/2}x^{1/2})+\O_p(n^{-1/6}x)+\O_p(n^{-1/3}),
\]
where the $\O_p$ terms do not depend on $x$.
Similar to~\eqref{eq:lwbound}, one obtains
\[
\widetilde{\mathbbm{Z}}_n^{\lambda}(x)-M_1x
\geq
\frac18\lambda_0'(x_0)x^2-S_n^2+o_p(1),
\]
for $M_2\leq x\leq n^{1/3}M_0$, where the $o_p$-term does not depend on $x$,
which yields
\[
\mathbbm{P}
\left(
\widetilde{\mathbbm{Z}}_n^{\lambda}(x)-M_1x\leq 0,
\text{ for some } M_2\leq x\leq n^{1/3}M_0\right)
\leq \varepsilon.
\]
In the case $x>n^{1/3}M_0$, similar to~\eqref{eq:convexity MLE},
Theorem~\ref{theorem:conv_breslow} and~\eqref{eq:marshall's lemma} yield
\begin{equation}
\label{eq:convexity LS}
\begin{split}
\BR(x_0+n^{-1/3}x)-\BR(x_0+M_0/2)
&\geq
\GCM(x_0+n^{-1/3}x)-\GCM(x_0+M_0/2)+o(1)\\
&\geq
\ho(x_0+M_0/2)(n^{-1/3}x-M_0/2)+o(1).
\end{split}
\end{equation}
This leads to
\[
\widetilde{\mathbbm{Z}}_n^{\lambda}(x)-M_1x
\geq
\widetilde{\mathbbm{Z}}_n^{\lambda}(n^{1/3}M_0/2)-M_1n^{1/3}M_0/2,
\]
from which we conclude
\[
\mathbbm{P}
\left(
\widetilde{\mathbbm{Z}}_n^{\lambda}(x)-M_1x\leq 0,
\text{ for some } x>n^{1/3}M_0\right)
\leq \varepsilon.
\]
This completes one part of the proof of~\eqref{eq:bounded Utilde lambda}
for a nondecreasing $\lambda_0$.
The other part is shown similarly.

For~\eqref{eq:bounded Utilde f}, using that $\widetilde{U}_n^{f}$ is nonincreasing,
similar to~\eqref{eq:Un to Zn MLE}, we first have
\[
\mathbbm{P}\left( n^{1/3}
\left[
\widetilde{U}_n^{f}( f_0(x_0)+n^{-1/3}M_1)-x_0\right]>M_2 \right)
\leq
\mathbbm{P}\left( \widetilde{\mathbbm{Z}}_n^{f}(x)-M_1x \geq 0, \text{ for some } x\geq M_2 \right),
\]
Next, according to~\eqref{eq:Znf in Znlambda expansion},
\eqref{eq:bound diff Zn LS ML} and~\eqref{eq:bound diff Wnhat x},
we obtain
\[
\begin{split}
\widetilde{\mathbbm{Z}}_n^{f}(x)-M_1x
&=
-
(1-F_0(x_0))n^{2/3}(\mathbb{P}_n-P)g(\cdot,n^{-1/3}x)\\
&\qquad-
(1-F_0(x_0))n^{2/3}Pg(\cdot,n^{-1/3}x)
-
\frac12(1-F_0(x_0))\lambda_0(x_0)^2x^2-M_1x\\
&\qquad+
\O_p(n^{-1/3})+\O_p(n^{-1/2}x^{1/2})+o_p(x)
+o_p(x^2),
\end{split}
\]
where the $\O_p$ and $o_p$ terms do not depend on $x$
and where $Pg(\cdot,n^{-1/3}x)$ is given in~\eqref{eq:Pg expansion}.
Now, choose $\delta=-f'(x_0)/(8(1-F_0(x_0)))>0$ in~\eqref{eq:Lemma41 Kim Pollard}, so that according to
Lemma~4.1 in~\cite{APPENDIXkimpollard:1990},
\[
-(1-F_0(x_0))n^{2/3}(\mathbb{P}_n-P)g(\cdot,n^{-1/3}x)
\leq
-\frac18 f_0'(x_0)x^2+S_n^2,
\]
for $n^{-1/3}x\leq M_0$ and $S_n^2=\O_p(1)$.
Furthermore, from~\eqref{eq:Pg expansion} together with~\eqref{eq:relation lambda' and f'}, it follows that we can find a $K>0$ such that for any $x>K$,
\[
-
(1-F_0(x_0))n^{2/3}Pg(\cdot,n^{-1/3}x)
-
\frac12(1-F_0(x_0))\lambda_0(x_0)^2x^2-M_1x
<
\frac14 f_0'(x_0)x^2,
\]
for $n$ sufficiently large.
Similar to~\eqref{eq:lwbound} we have for $x\geq M_2\geq K$,
\[
\widetilde{\mathbbm{Z}}_n^{f}(x)-M_1x
\leq
\left(
\frac18f_0'(x_0)+o_p(1)
\right)x^2
+
S_n^2 +o_p(1),
\]
where the $o_p$ terms do not depend on $x$, which leads to
\[
\mathbbm{P}\left(
\widetilde{\mathbbm{Z}}_n^{f}(x)-M_1x \geq 0,
\text{ for some } M_2\leq x\leq n^{1/3}M_0\right)
\leq
\varepsilon,
\]
for $n$ sufficiently large.
In the case $x>n^{1/3}M_0$, first, similar to~\eqref{eq:marshall's lemma}, we can obtain
that for any $0<M<\tau_H$,
\[
\sup_{x\in[0,M]}\left| \widetilde{F}_n(x)-F_0(x) \right|
\leq
\sup_{x\in[0,M]}\left|F_n(x)-\Lambda_0(x) \right|,
\]
which then similar to~\eqref{eq:convexity LS} together with Corollary~\ref{cor:consistency_F_n} yields
\begin{equation}
\label{eq:convexity dens}
\begin{split}
F_n(x_0+n^{-1/3}x)-F_n(x_0+M_0/2)
&\leq
\widetilde{F}_n(x_0+n^{-1/3}x)-\widetilde{F}_n(x_0+M_0/2)+o(1)\\
&\leq
\fo(x_0+M_0/2)(n^{-1/3}x-M_0/2)+o(1).
\end{split}
\end{equation}
This leads to
\[
\widetilde{\mathbbm{Z}}_n^{f}(x)-M_1x
\leq
\widetilde{\mathbbm{Z}}_n^{f}(n^{1/3}M_0/2)-M_1n^{1/3}M_0/2,
\]
from which we conclude
\[
\mathbbm{P}
\left(
\widetilde{\mathbbm{Z}}_n^{\lambda}(x)-M_1x
\geq 0,
\text{ for some } x>n^{1/3}M_0\right)
\leq \varepsilon.
\]
This completes one part of the proof of~\eqref{eq:bounded Utilde f}.
The other part is shown similarly.

Finally, the proof of~\eqref{eq:bounded Uhat lambda} and~\eqref{eq:bounded Utilde lambda} in the case of a nonincreasing $\lambda_0$
is along the lines of the proof of~\eqref{eq:bounded Utilde f}, combined with arguments used for the proof
of~\eqref{eq:bounded Uhat lambda} and~\eqref{eq:bounded Utilde lambda} in the nondecreasing case.
\end{proof}

%\bibliographystyle{acm}
%\bibliography{cox_RL_TN}
%\nocite{groeneboomwellner}
%\nocite{huangwellner:1995}
%\nocite{kimpollard:1990}

%\bibliographystyle{acm}
%\bibliography{cox_RL_TN}

\end{document}